\documentclass[onecolumn,preprintnumbers,amsmath,amssymb]{revtex4}
\usepackage{graphicx}%  Include   figure  files  \usepackage{dcolumn}
\usepackage{bm}
\usepackage{color}
\begin{document}

%\preprint{APS/123-QED}

\title{On a dynamical approach to some prime number sequences}
% Force line breaks with \\

%\author{B Luque, L Lacasa}
%\email{lucas@dmae.upm.es}
%\affiliation{Dpto. de Matem\'atica Aplicada y Estad\'istica\\
%ETSI Aeron\'auticos\\ Universidad Polit\'ecnica de Madrid.}%

\author{Lucas Lacasa$^1$, Bartolo Luque$^2$, Ignacio G\'{o}mez$^2$, and Octavio Miramontes$^{2,3}$}
\affiliation{
$^1$School of Mathematical Sciences, Queen Mary University of London, Mile End E14NS London, (United Kingdom)\\
$^2$Escuela Tecnica Superior de Ingenier\'{i}a Aeron\'{a}utica y del Espacio,  Universidad Polit\'ecnica de  Madrid,   Madrid  28040 (Spain)\\
$^3$Departamento  de   Sistemas
Complejos, Instituto  de F\'isica, Universidad  Nacional Aut\'onoma de
M\'exico, 04510 DF (Mexico)
}%

\date{\today}
%  It is  always \today,  today, %  but any  date  may be
%  explicitly specified

\begin{abstract}
In this paper we show how the cross-disciplinary transfer of techniques from Dynamical Systems Theory to Number Theory can be a fruitful avenue for research. We illustrate this idea by exploring from a nonlinear and symbolic dynamics viewpoint certain patterns emerging in some residue sequences generated from the prime number sequence. We show that the sequence formed by the residues of the primes modulo $k$ are maximally chaotic and, while lacking forbidden patterns, display a non-trivial spectrum of Renyi entropies which suggest that every block of size $m>1$, while admissible, occurs with different probability. This non-uniform distribution of blocks for $m>1$ contrasts Dirichlet's theorem that guarantees equiprobability for $m=1$. We then explore in a similar fashion the sequence of prime gap residues. This sequence is again chaotic (positivity of Kolmogorov-Sinai entropy), however chaos is weaker as we find forbidden patterns for every block of size $m>1$. We relate the onset of these forbidden patterns with the divisibility properties of integers, and estimate the densities of gap block residues via Hardy-Littlewood $k$-tuple conjecture. We use this estimation to argue that the amount of admissible blocks is non-uniformly distributed, what supports the fact that the spectrum of Renyi entropies is again non-trivial in this case. We complete our analysis by applying the Chaos Game to these symbolic sequences, and comparing the IFS attractors found for the experimental sequences with appropriate null models.
\end{abstract}

%\pacs{05.70.Fh,  05.10.Ln,  02.10.De,  64.60.-i,  64.60.Cn,  64.60.Fr, 89.20.Ff}

\maketitle %%%

\section{Introduction}
Number Theory is a millennial branch of pure mathematics devoted to the study of the integers, whose implications and tentacles -despite the expectations of theorists like Dickson or Hardy- not only pervade today almost all areas of mathematics but are also at the basis of many technological applications. A particularly interesting and fruitful bridge to travel is the one linking dynamical systems concepts with number-theoretic ideas, and indeed several well established subfields in pure mathematics lie at the interface between number theory and dynamics, namely arithmetic dynamics, dynamics over finite fields or Lie groups, symbolic dynamics, to cite some \cite{Schroeder, 5}.\\

\noindent Quite distant from these in style and focus, physicists have also history of dealing with dynamics, originating in Classical Mechanics and more recently encompassing areas such as Chaos Theory or Complexity Science. As a matter of fact, in the last decades physicists with their tools and experimental inclination have started to look at number theoretic sequences as experimental measurements extracted from some hidden underlying dynamics. Moreover, from the celebrated coincidence in 1972 between Montgomery's work on the statistics of the spacings between Riemann zeta zeros and Dyson's analogous work on eigenvalues of random matrices in Nuclear Physics \cite{Dyson, Montgomery, Odlyzko}, we have seen, somewhat unexpectedly, how number theory and physics have built bridges between each other. These connections range from the reinterpretation of the Riemann zeta function as a partition function \cite{1} or the focus of the Riemann hypothesis via quantum chaos \cite{2}, to multifractality in the distribution of primes \cite{3}, computational phase transitions and criticality originating in combinatoric problems \cite{4, lucasPRE, lucasPRL, lucasPRE2, lucasPRE3}, or stochastic network models of primes and composites \cite{networks} to cite only a few examples (see \cite{5} for an extensive bibliography).\\

\noindent In this work we aim to illustrate this fertile cross-disciplinary transfer of ideas and tools by tackling from a dynamical point of view some important sequences that emanate from the prime number sequence, with the aims of exploring its underlying structure. By using tools originally devised to describe turbulent fluids or to generate fractal patterns, we find that these sequences show compelling signs of chaotic behavior while hiding some unexpected structure. In an effort to elucidate this dynamical interpretation, we then link these results with the number theoretical properties of these sequences.\\

\begin{figure}[ht!]
\centering
\includegraphics[scale=1.57]{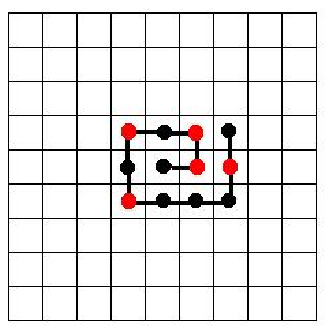}
\includegraphics[scale=0.393]{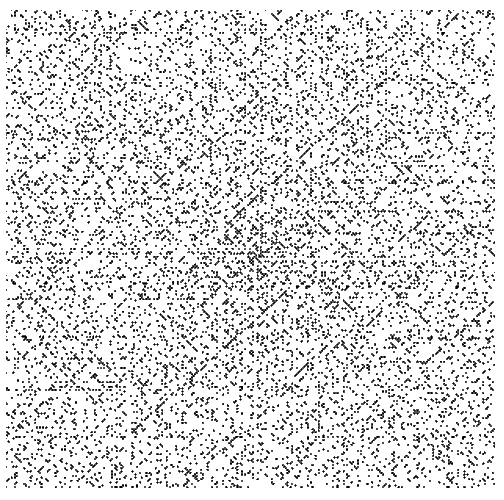}

\caption{(left)The Ulam spiral for the sequence 1, 2 $\dots$ 12, with primes in red. (right) A full 200 x 200 Ulam spiral showing the primes as individual black dots. Notice the emergent diagonal pattern where the primes tend to accumulate.}
\label{fig:ulam_spiral}
\end{figure}

\noindent {\bf Prime spirals and the residue classes of linear congruences. }
The Ulam spiral \cite{gardner1964remarkable} is make by writing integers in a square spiral and marking the particular position of the prime numbers (see figure \ref{fig:ulam_spiral} for an illustration). Using this representation, S. Ulam found that prime numbers tend to distribute and appear mostly on the diagonals of the spiral.
Commonly the primes are marked in a different color like in figure \ref{fig:ulam_spiral} (left) in order to distinguish them. The diagonal-like accumulation patterns of primes can be explained by noticing that these correspond to solutions of prime-generating polynomial
equations such as the well known Eulerian form
 $f(n) = 4n^2 + bn + c. $\\
 %{\color{red} Give more details here, maybe talk about Hardy and Littlewood's Conjecture F.}
%Hanh and  Sachs generalized the representation of the Ulam spiral in the so-called Sachs spirals \cite{hahn2008distribution}. Furthermore, there is an infinite number of these geometric sequences as shown recently \cite{tao2008primes}.\\

\begin{figure}[ht!]
\centering
\includegraphics[scale=1.65]{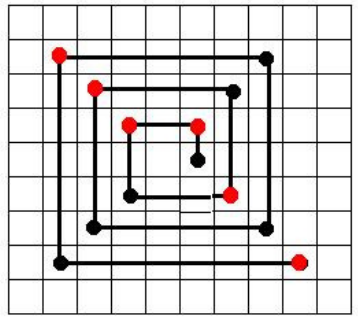}
\includegraphics[scale=0.45]{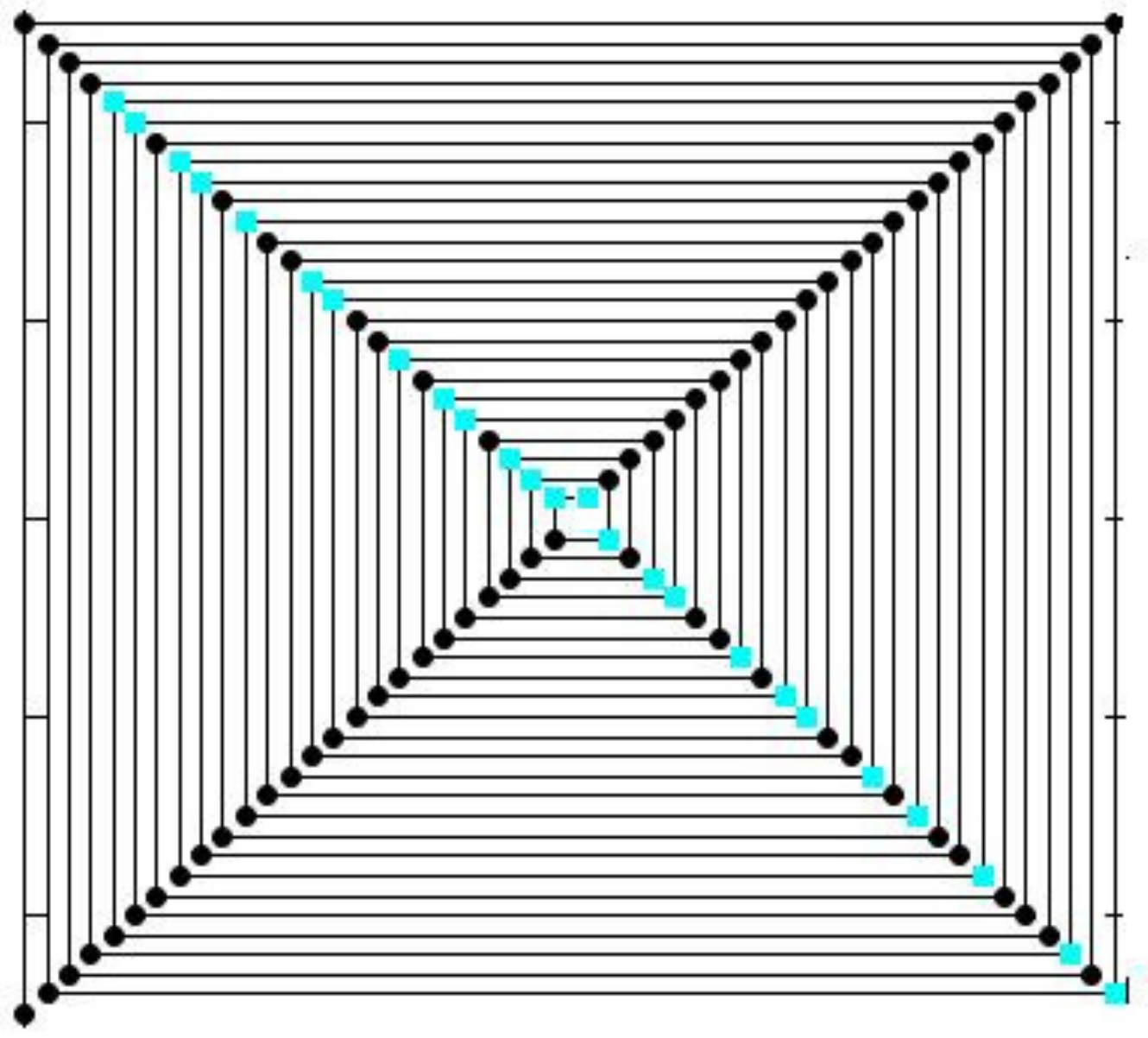}
\caption{(Left panel) The diagonal spiral for the sequence 1, 2 $\dots$ 13, with primes in red. (Right panel) The diagonal spiral for $N=\{2 \dots 100\}$, with primes in blue. The upper-left spiral arm agglutinates primes of the form $4n+3$, called Gaussian primes. The bottom-right agglutinates the rest of the primes with form $4n+1$, called Pythagorean primes.}
\label{fig:diag_spiral_detail}
\end{figure}

\noindent A simpler spiral that can be constructed to represent the sequence of natural numbers is the so-called diagonal spiral, where the spiral rotates off from an initial point and the interval of the distance between successive numbers is increased by one each two intervals as shown in Figure \ref{fig:diag_spiral_detail}.
The resulting spiral accommodates all integer in four diagonals as shown in Figure \ref{fig:diag_spiral_detail} (right panel), for $N=100$. With the exception of 2, all prime numbers are accommodated on two opposite diagonals due to the fact that all primes apart besides 2 are odd numbers. As this spiral is a simple geometric illustration of performing a modular operation (modulo 4), then each arm of the spiral agglutinates integers of the form $4n+ b$, with $n$ integer and $b=0,1,2,3$. The upper-left spiral arm agglutinates primes such as $\{3, 7, 11, 19, 23, 31,...\}$ of the form $4n + 3$, i.e. congruent to 3 modulo 4, and are called {\it Gaussian primes} (this sequence is catalogued as A002145 in OEIS). On the other hand, the bottom-right diagonal agglutinates the rest of the primes, e.g. $\{5, 13, 17, 29, 37, 41..\}$ of the form $4n+1$, i.e. congruent to 1 module 4, and are called {\it Pythagorean primes} (sequence A002144 in OEIS). Besides the prime $2$, all primes can be expressed as $4n+1$ or $4n+3$ for some $n\geq 0$. Our analysis will start by considering a symbolic sequence $p(n) \mod 4$ where $p(n)$ is the n-th prime, with two symbols $\{A,B\}$ such that $A\equiv p(n) \mod 4=1$, $B\equiv p(n) \mod 4=3$ \cite{Hao}. By virtue of Dirichlet's theorem these two symbols occur infinitely often along the prime number sequence.\\

\noindent Incidentally, note that this initial mapping can be performed with arbitrarily different modulus $k$; i.e. $kn+b$ with $b=0,1,2,...,k-1$. Suppose that we construct the residue classes formed by taking the primes modulo $k$. Dirichlet's theorem on arithmetic progressions proves that all but a finite number of primes are congruent with $\phi( k)$ different residue classes, where $\phi(\cdot)$ is the Euler totient function \cite{Schroeder}. For instance, since $\phi(3)=\phi(4)=\phi(6)=2$, this means that for $k=3,4$ and $6$ the symbolic sequence constructed from taking the residue classes of primes modulo $k$ is essentially formed by just two symbols (note however that there is always a finite transient, for instance for $k=4$ the first prime $p(1)=2$ is the only prime which is not in the residue classes $1$ and $3$). In this sense, the sequence formed by Pythagorean and Gaussian primes is just a convenient particular case to distinguish both types of primes, but one can do this for arbitrary modulus. For instance, for modulus $k=5,8,10,12$ the totient function is $4$, so the resulting sequence will have four symbols in these cases, and so on.In this work we will only consider prime residue sequences with two symbols obtained via taking the primes modulo $k$, with $k=3$ and $6$ (Dirichlet's theorem guarantees that the density of both symbols is positive in every case).\\

\begin{figure}[ht!]
\centering
\includegraphics[scale=0.6]{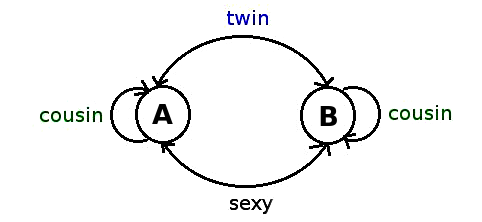}
\caption{A network representation of the prime numbers sequence and state transitions. Prime sequence can be seen as the result of a Markov Chain defined over a two-state network with certain transition rates. In figure \ref{fig:two_branches} we depict an illustration of such a Markov Chain process.}
\label{fig:two_branches_1}
\end{figure}

\begin{figure}[ht!]
\centering
\includegraphics[scale=0.6]{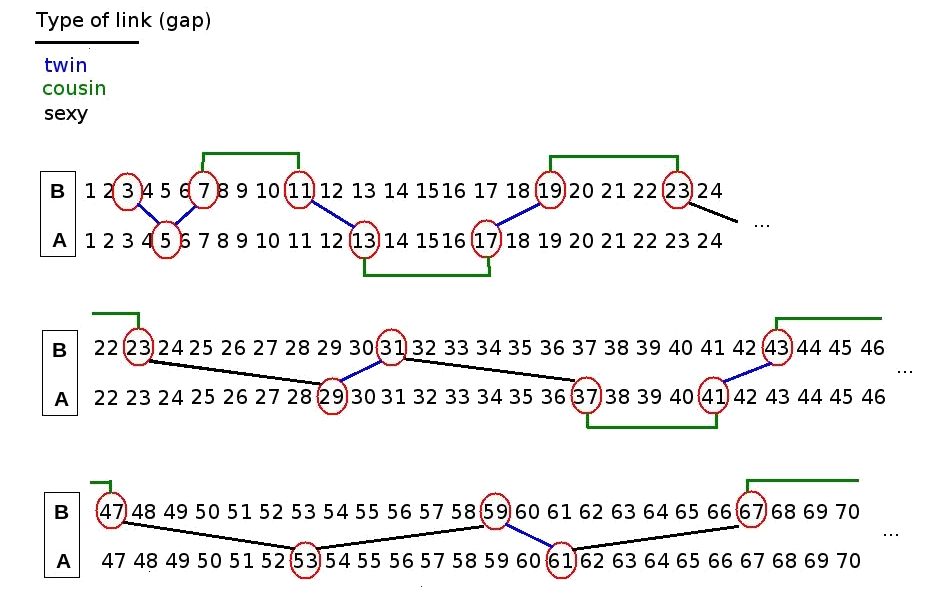}
\caption{Transitions between consecutive primes, labelled in terms of twins, cousins and sexies.}
\label{fig:two_branches}
\end{figure}

\noindent {\bf The transitions between Pythagorean and Gaussian primes.} Once the sequence of Gaussian and Pythagorean primes (modulus $k=4$) has been extracted, we can consider the sequence of transitions between these two classes. This new symbolic sequence now has 4 symbols, one per transition, i.e. the set $\{AA,AB,BA,BB\}$. We construct this sequence by sliding a block-2 window and assigning to each pair of consecutive symbols a new symbol, for instance the consecutive primes $3$ and $5$ map into the symbol $AB$, because $3$ is a Gaussian prime ($A$) and $5$ is a Pythagorean prime ($B$). This new sequence will be called the {\it transition sequence}. The reason for doing this is inspired in envisaging primes as the result of a Markov Chain defined over a two-state network with certain transition rates, which would in this case be equivalent to the frequencies of each block-2 string (see figure \ref{fig:two_branches_1}).\\

\noindent {\bf Residue classes of gaps mod $k$: twins, cousins and sexies.} Two consecutive primes on the integer line separated by a gap $g=2$ are known as {\it twin} primes (for example $3, 5$ are twins).
Consecutive primes on the integer line separated by a gap $g=4$ or $6$ are known as {\it cousin} primes and {\it sexy} primes respectively (for example $7, 11$ are cousins and $23,29$ are sexies). An illustration of the occurrence of twin, cousin and sexy prime pairs is shown in figure \ref{fig:two_branches}.
Gaps between primes, ever odd except for the number $2$ and $3$ the only gap with $g=1$, can grow arbitrarily large, so for the sake of studying the properties of gaps with the tools of dynamical systems theory we need to make the sequence stationary. Again, we do this by taking modulo $k$ and symbolize the gaps according to the residue class to which they correspond. For $k=6$, one can accordingly classify the gaps into three families: those congruent with $0 \mod 6$ will be called sexy-like, those congruent with $2 \mod 6$ will be called twin-like, and those congruent with $4 \mod 6$ will be called cousin-like. Incidentally, note that transitions between the Pythagorean and Gaussian branches are only due to twin-like and sexy-like gaps. The $k$-tuple conjecture by Hardy and Littlewood \cite{hardy, ares} indeed predicts that all admissible $k$-tuple of primes occur infinitely often, hence this conjecture suggests that the density of the residue classes $0,2$ and $4$ is positive accordingly. Our third type of sequence under analysis in this work will be the sequence of gap residues $g(n):= (p(n+1)-p(n)) \mod 6$.\\

\noindent The rest of the paper goes as follows: in section II we recall some important tools originating in dynamical systems and statistical physics for the description of disordered systems from sequence measurements. In section III we apply these methods to the three types of sequences described above: (i) the symbolic sequence of prime residues $p(n) \mod 4$ with symbols $A$ (Pythagorean primes) and $B$ (Gaussian primes), (ii) the transition sequence among them with four symbols: $AA$, $AB$, $BA$ and $BB$ and (iii) the gap sequence $g(n) \mod 6$ wit three symbols, $A$ for sexy primes, $B$ for twin-like primes, and $C$ for cousin-like primes. And interpret these results from both a dynamical and number theoretical viewpoint.
In section IV we discuss our findings and conclude.

\section{Tools from nonlinear and symbolic dynamics}

\subsection{Enumerating blocks: Spectrum of Renyi entropies}
\noindent In this subsection we show how to explore in a quantitative and systematic way the abundance of blocks of symbols that appear in each symbolic sequence, and its dynamical interpretation. A block of size $m$ constitutes a string of $m$ consecutive symbols in the sequence: is every possible block appearing in the sequence? If so, with what frequency?\\

\noindent Let us consider a dynamical system $x_{t+1}=F(x_t)$, $F:{\cal X}\subset  \mathbb{R}\to {\cal X} \subset \mathbb{R}$ and consider a symbolic sequence ${\cal S}=(s_0,s_1,\dots)$ extracted from a trajectory $x_0,F(x_0),F^2(x_0),\dots$ of this map via a certain partition of the phase space ${\cal P}$ into $p$ symbols. In this infinite sequence we are interested in the statistical properties of a generic block of $m$ consecutive symbols ${\bf s}=[s_1\dots s_m]$. To fully describe the statistics of these blocks we shall consider the whole spectrum of so-called Renyi dynamical entropies $h(\beta)$, with $\beta \in \mathbb{R}$ \cite{beck} (Renyi entropies, for short) that weight the measure, the frequency of each block in different ways. There are two particular cases that, given their paramount importance, should be highlighted before a general formula is introduced.\\

\noindent For $\beta=0$, $h(0)$ is called the topological entropy. If ${\cal A}(m)$ is the set of all admissible blocks (present in the dynamics) of length $m$, then the topological entropy can be defined in terms of $|{\cal A}(m)|$ as:
\begin{equation}
|{\cal A}(m)|\sim e^{m h(0)}
\label{htop}
\end{equation}
In other words, $h(0)$ describes how many new different admissible blocks can appear in the symbolic sequence as we increase its length $m$. We can compute $h(0)$ as:
\begin{equation}
h(0)=\sup_{\cal P}\lim_{m\to \infty} \frac{1}{m}\log{|{\cal A}(m)|}.
\label{h0}
\end{equation}

\noindent There is no metric underlying this quantity, we are only counting admissible blocks, hence the label {\it topological}. This is the symbolic analog of the rate of new trajectories that are admissible in a dynamical system as time increases.\\

\noindent For $\beta=1$, $h(1)$ reduces to the symbolic version of the well known Kolmogorov-Sinai entropy:
\begin{equation}
h(1)=\sup_{\cal P}\lim_{m\to \infty} \frac{-1}{m}\sum_{{\cal A}(m)} P({\bf s})\log P({\bf s}).
\label{ks}
\end{equation}
Where the summation goes through all the admissible configurations {\bf s} and P({\bf s}) is the probability of each one of them. This entropy weights logarithmically the measure of each type of block, it is thus a metric entropy.
A traditional definition of a chaotic process is associated with a finite, positive value of $h(1)$, and under suitable conditions this quantity is, according to Pesin identity \cite{beck}, equivalent to the sum of positive Lyapunov exponents of the underlying (hidden) chaotic system.\\

\noindent For $\beta>1$, $h(\beta)$ defines a spectrum of Renyi entropies:
\begin{equation}
h(\beta>1)=\sup_{\cal P}\lim_{m \to \infty} {\frac{1}{1-\beta}\frac{1}{m} \log \sum_{{\cal A}(m)}}P({\bf s})^\beta.
\label{renyi}
\end{equation}
It is known that $h(\beta)$ is a monotonically decreasing function of $\beta$, such that $h(\beta+1)\leq h(\beta)$. For some chaotic processes such as the binary shift map \cite{Schroeder, Hao}, one can show that $h(\beta)$ is independent of $\beta$, one then says that the spectrum collapses to a single value. If $h(\beta)$ indeed depends on $\beta$ we say that the system has a non-trivial spectrum of Renyi entropies.\\

\noindent Two comments are in order. First, note that $h(\beta)$ are strictly speaking not simple entropies but {\it entropy rates}, that is, they are intensive quantities in the statistical mechanics sense. For convenience, we will also use the following notation
$$h(\beta)=\lim_{m\to \infty} \frac{H_m(\beta)}{m},$$
where $H_m(\beta)$ are the Renyi block entropies, also called partial entropies. These are indeed extensive quantities, and the entropy rate is taken by normalizing these over the block size. Second, note that all these entropies in general depend on the partition ${\cal P}$ performed in the phase space of the underlying dynamical system. If the partition is performed homogeneously, then this essentially depends on the number of symbols, but this is just a particular subset of all possible partitions. Strictly speaking, in order to find the correct value of $h(\beta)$ one needs to take the supremum over all possible partitions. We will obviate this step, relying on the fact that for the so-called generation partitions, these quantities already reach their maximum and therefore the supremum does not need to be taken.

\subsection{IFS and the Chaos Game}
\begin{figure}[ht!]
\centering
\includegraphics[scale=0.27]{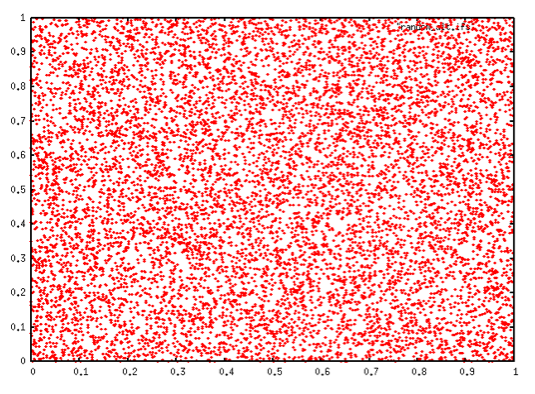}
\includegraphics[scale=0.2]{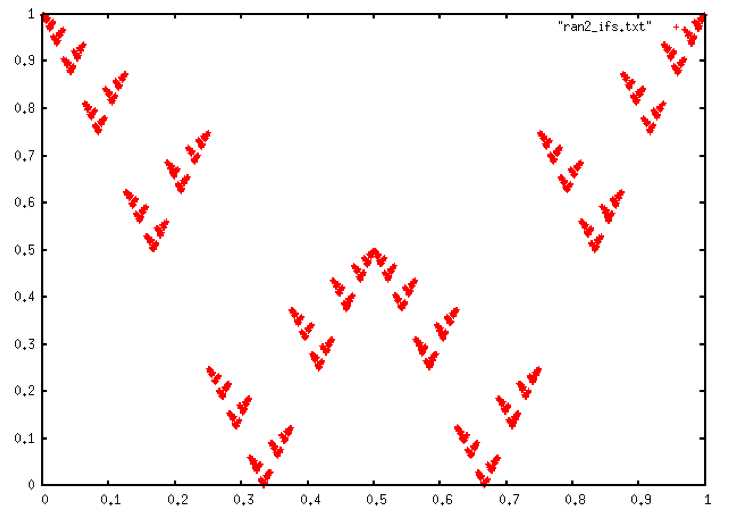}
\includegraphics[scale=0.12]{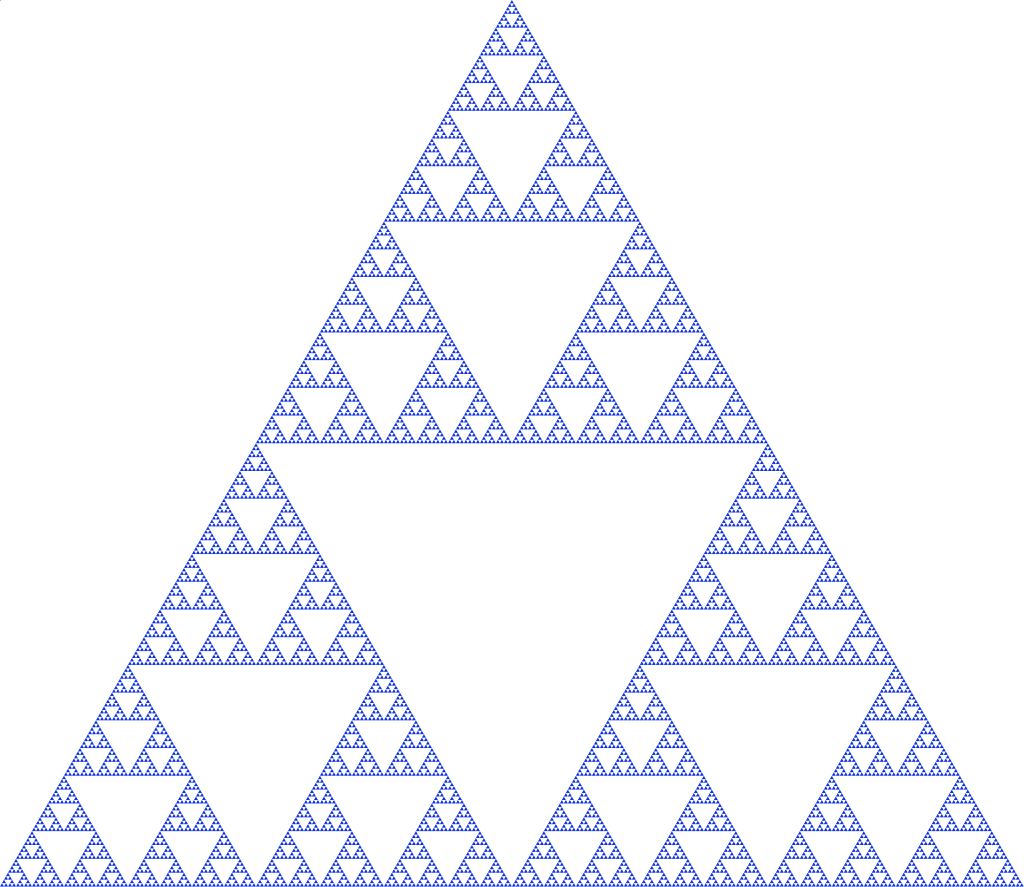}
\caption{IFS attractor for (left panel) a type I null model with $p=4$ symbols, (middle panel) a type II null model with $p=4$ symbols, (right panel) a type I null model with $p=3$ symbols.}
\label{fig:nullmodeltransition}
\end{figure}

The so-called Chaos Game \cite{barnsley,peitgen,jeffrey} is a simple iterative method to construct fractals. Originally, one starts with an initial point $x_0$ inside a $p$-vertex polygon and builds a new point $x_1$ by selecting {\it at random} one of the $p$ vertices of the polygon, drawing the segment connecting $x_0$ and the vertex, and finding the location of the point in the segment whose distance to the vertex is a certain factor of the distance between $x_0$ and the vertex. This transformation is  an affine transformation, whose iteration defines a so-called Iterated Function System (IFS), a well-known iterative method to create fractals \cite{barnsley}.
For $p=3$ (regular triangle) and a factor $1/2$, the attractor of this simple iterative process is the celebrated Sierpinski gasket, whereas for $p=4$ (square) and a factor $1/2$, the trajectory is space filling and the attractor is the whole square (see figure \ref{fig:nullmodeltransition}).\\

\noindent Now, suppose that the process of picking new vertices is not at random, but follows a certain pattern. In that case the attractor of the IFS is typically a subset of the attractor found by making the Chaos Game. For instance, it is well known that for $p=4$ and a factor $1/2$, the attractor of the Chaos Game is the whole square, but if restrictions on what vertices can be selected based on previous history are applied; i.e. if we introduce temporal correlations in the stochastic process of vertex selection, in several cases the new attractor is a fractal subset of the square. Given a symbolic sequence, where symbols can take $p$ different values, one can therefore apply the IFS defined above where in each iteration the selected vertex is not chosen at random but is given by the next element in the symbolic sequence. Non-random sequences will then in general yield an IFS attractor different from the one obtained applying the Chaos Game, and the concrete shape would be given by the precise way we have to prune a random series to obtain the symbolic sequence we are working with.\\

\noindent This argument provides a way to construct a geometric representation of the patterns underlying non-random symbolic sequences, which can be disentangled by comparing them with the object retrieved by applying the Chaos Game (that would be called the null model, see next section for details). In what follows we will explore the geometric shapes we find by applying the Chaos Game not on random sequences but following the specific sequences extracted from the primes.

%any real-valued one dimensional sequence $X$ can be embedded --and its underlying dynamic exhibited-- into a two-dimensional space bounded inside a $p \geq 3$ vertex polygon, a procedure known commonly as the Chaos Game \cite{barnsley,peitgen,jeffrey}. This procedure implies symbolising $X$ into a discrete $p$-valued symbolic sequence $S$ and then perform an iterated matrix affine transformation on a vector plus a translation (so-called Iterated Function System procedure, or IFS).
%The IFS analysis performed on the gap residue sequence coming from eq.\ref{eq:mod6} is shown in figure \ref{fig:IFSgaps}, where the result of a type I null model is plotted as well. We find that the prime gaps selects a subset of the attractor generated by the null model. This particular subset could be put in correspondence with the subset of admissible sequences found above, and in some sense the IFS attractor plays the geometric role of the distribution of forbidden patterns.

\section{Results}
\subsection{Null models: what to expect if primes lack structure?}
What should we expect from a null model? Here we are interested in evaluating whether our prime sequences hides some degree of correlations or, conversely, these sequences cannot be distinguished from a totally random process. As such, our null model should consist of a random uncorrelated sequence. However, we have several possibilities to do that.

\subsubsection{Type I: Symbols i.i.d. with uniform and non-uniform probability densities}
In such a model, $s$ is a random variable extracted at each time independently with a certain probability distribution from $\{s_1,\dots,s_p\}$, where $p$ is the number of symbols. For instance, in our case one could draw symbols {\it uniformly} from either $\{A,B\}$ or $\{AA,AB,BA,BB\}$. That is, each of the two or four symbols having a probability $1/2$ or $1/4$ respectively, to generate a random sequence, which would act as a null model for both the prime residue sequence and the transition sequence respectively.
It is easy to prove that in this case, $H_m(\beta)= \log p$, independently of $\beta$ and $m$, something that was checked numerically for $\beta=1$ in figure \ref{fig:h1}. This means that in the this kind of null model, all blocks of size $m$ are possible and appear in the symbolic sequence with equal probability. The case $\beta=0$ is trivial to prove, here we give a sketch of the proof for $h(1)$: in the case all blocks of size $m$ are equiprobable, we have $P({\bf s})=1/|{\cal A}(m)|$ and then
$$H(1)=\sum_{{\cal A}(m)} P({\bf s})\log P({\bf s})=\log |{\cal A}(m)|.$$
In the case that all possible blocks of size $m$ are admissible, the set of admissible sequences is equal to the set of all possible sequences ${\cal N}(m)$. We thus have ${|\cal A}(m)|={|\cal N}(m)|=p^m$ hence
$$h(1)=\lim_{m\to \infty} \frac{1}{m}\log(p^m)=\log p.$$
Furthermore, we can also define a similar stochastic process where now we draw symbols at random from the set of symbols not in a uniform way, but instead with a probability that indeed matches the empirical one. This would be a null model that respects marginal distributions.
For instance, for primes residues where the symbol set is $\{A,B\}$, while asymptotically $P(A)=P(B)=1/2$, for finite sequences estimations are systematically larger for $A$, something known as the Tchebytchev bias. For the transition sequence where the symbol set is $\{AA,AB,BA,BB\}$, empirically we find $P(AB)=P(BA)\approx 0.30, \ P(AA)\approx 0.21,\ P(BB)\approx 0.19$: this might account for this bias as $P(A)= P(AB)+P(AA) > P(BA) + P(BB)= P(B)$ (note that for convenience we have not made explicit in the notation the difference between asymptotic and estimated probabilities).\\

In these cases, the null model accounts for this non uniformity in the abundance of symbols and generates a sequence with non-uniform symbol frequencies which, for finite $m$, has by construction metric entropies which are smaller or equal to the uniform counterpart (this holds as a direct consequence of the second Khinchin axiom). Consider the case $\beta=1$. For blocks of size $m=1$ the entropy $H_1(1)$ is:
$$H_1(1)=-\sum_{{\cal A}(1)} P(s_1)\log P(s_1)\approx 1.366 < \log 4 \approx 1.386.$$
This value holds constant $\forall m$ and one can very easily prove that in this case $h(1)=\lim_{m\to \infty} H_m(1)/m=H_1(1)$. A sketch of the proof is as follows: consider for simplicity the case $m=2$, in that case:
$$H_2(1)=\sum_{{\cal A}(2)} P({\bf s})\log P({\bf s})\equiv \sum_{s_1=1}^p \sum_{s_2=1}^p P(s_1,s_2)\log P(s_1,s_2)$$
Since the sequence is by definition uncorrelated we can factorize $P(s_1,s_2)=P(s_1)P(s_2)$ and thus
$$H_2(1)=\sum_{s_1=1}^p \sum_{s_2=1}^p P(s_1,s_2)\log P(s_1,s_2)=\sum_{s_1=1}^p P(s_1)\log P(s_1) \sum_{s_2=1}^p P(s_2) + \sum_{s_1=1}^p P(s_1)\sum_{s_2=1}^p P(s_2)\log P(s_2)= 2H_1(1)$$
A simple induction concludes the proof. Accordingly, the null model for the prime configuration depicted above yields $H_m/m \approx 1.366 \ \forall m$, and in general $h(\beta)=H_1(\beta)\leq \log 4$ (the inequality saturates for $\beta=0$ as in this case no metric aspects are concerned, whereas for other values of $\beta$ we expect the upper bound to be reasonably tight).\\

\noindent Now, what type of IFS attractor is produced by a type I null model? If one performs the Chaos Game for a square (4 vertices) and a contraction factor $1/2$, the process is space-filling, i.e. the attractor is $[0,1]^2$ (the attractor is the interior of a square with all points visited with equal probability). Hence performing the Chaos Game on a type I null model with $p=4$ symbols is a space-filling IFS (left panel of figure \ref{fig:nullmodeltransition}). On the other hand, for $p=3$ a type I null model yields the Sierpinski gasket as the attractor of the IFS (right panel of figure \ref{fig:nullmodeltransition}), point out no correlations.

\subsubsection{Type II: transitions}
After a careful thought it is straightforward to see that the preceding null models are reasonable to explore the lack of correlations in sequences such as the primes or gaps residues, but not for the prime transition sequence, made out of $\{AA,AB,BA,BB\}$. The reason is simple: suppose that in the sequence at some point we find the symbol $AB$. The subsequent symbol will describe the transition from $B$ to the next prime, and therefore this next symbol will necessarily start with the letter $B$ (that is, $BA$ or $BB$), in other words it will never be a symbol starting with the letter $A$. This is simply related to the way of constructing the transition sequence (by sliding a window of size $2$ with overlap) hence this process builds spurious correlations and forbidden transitions. As an example: the following transitions are forbidden:\\
$\{ AA\to BA, \ AA \to BB,\  AB \to AA,\  AB \to AB,\  BA \to BA,\  BA \to BB,\  BB \to AA, \ BB \to AB\}$.\\

\noindent Accordingly, plenty of possible blocks of $m$ symbols are actually forbidden simply by the way the sequence is created, not necessarily because there is an intrinsic correlation in the prime sequence. In order to find out whether this is the case, the correct null model to explore is as follows:
\begin{itemize}
\item Generate a random sequence via a type I null model with two symbols $A$ and $B$ (where asymptotically $p(A)=p(B)$, however for finite size series one should account for Tchebytchev bias),
\item Then construct the resulting sequence with $p=4$ symbols $\{AA,AB,BA,BB\}$.
\end{itemize}

\noindent It is easy to enumerate in this null model the number of possible $|\cal N|$, admissible $|\cal A|$ and forbidden $|\cal F|$ blocks of size $m$, with $|{\cal N}(m)|=|{\cal A}(m)|+|{\cal F}(m)|$, throught the iterations:\\
$$|{\cal F}(m)|= 4|{\cal F}(m-1)| + 2|{\cal A}(m-1)|, \ |{\cal A}(m)|=|{\cal N}(m)|-|{\cal F}(m)|.$$
With $|{\cal N}(m)|=4^m$ and $|{\cal F}(1)|=0$ we can solve the preceding equations and find:
\begin{equation}
|{\cal F}(m)|=2^m(2^m-2),\ |{\cal A}(m)|=2^{m+1}
\label{nullII}
\end{equation}
Accordingly, in this null models the number of forbidden patterns grows exponentially fast (faster than the number of admissible ones). The topological entropy of this model is analytically solvable:\\
$$h(0)=\lim_{m\to \infty}\frac{1}{m} \log {\cal A}(m)=\lim_{m\to \infty}{\frac{m+1}{m}\log 2}=\log 2,$$
and $H_m(0)=\log 2 + (\log 2)/m$. We find that this quantity converges to $\log 2$ after a monotonic decrease.\\

\noindent The IFS attractor of a type II ($p=4$) is plotted in the middle panel of figure \ref{fig:nullmodeltransition}. The attractor is not anymore the interior of the square (i.e. an attractor of dimension 2), but a fractal of lower dimension. This resulting fractal shape is the geometrical equivalent of systematically pruning from ${\cal N}(m)$ the spurious forbidden patterns discussed above.

\subsubsection{Type III: Cramer null model}
There is a third null model than one could. This is the classical Cramer random model \cite{Cramer}, which starts by creating a stochastic version of the prime number sequence by assigning to each integer $x$ a probability of being prime equal to $1/\log x$. This model is inspired after the prime number theorem, which states that the amount of primes under N is asymptotic with $N/\log N$. This latter quantity can be therefore understood as the expected number of primes under $x$. If the probability of a number only depends on its ordinal in the line, then the average number of primes under $N$ is
$$\pi(N)= \int_2^N \frac{dx}{\log x}=\text{Li}(N)\sim N/\log N$$
where $\textrm{Li}(x)$ is the offset logarithm integral. This model generates a list of pseudo-primes which statistically conform to the prime number theorem albeit being essentially uncorrelated from each other.
%Cramer assumed that, given $x$ the fact whether each integer in (1,x) is prime or composite follows from a Bernoulli trial, and then deduced that such trial happens with probability $1/\log x$, this way the estimated number of primes in a given interval (1,x) coincides with the formula given by the prime number theorem, but the set of pseudo-primes generated by this process lack any correlation structure.
Note that in this form, Cramer model is not applicable in our case, as it is very likely that such model would generate both even and odd pseudo-primes, therefore not congruent with the correct residue classes. We therefore slightly modify this model and assume that each even integer $x$ is a pseudo-prime with null probability, whereas each odd integer $x$ is a pseudo-prime with probability $2/\log x$. From the resulting list of pseudoprimes one can reconstruct the pseudo-prime residues, pseudo-prime transitions and pseudo-prime gap residues sequences.\\

\noindent We won't give more details at this point. It is just important to highlight that a deterministic process whose statistics are equivalent to the corresponding null model is totally indistinguishable from a purely random, uncorrelated process. The main aim of this work is to explore this hypothesis.

\subsubsection{Some fully chaotic maps}
\noindent For the sake of comparison with other deterministic processes, we close this section on null models by considering several fully chaotic maps: the logistic map $x_{t+1}=4x_t(1-x_t)$, the binary shift map $x_{t+1}=2x_t \mod 1$ and the tent map $x_{t+1}=2 \min\{x_t,1-x_t\}$ where in every case $x \in [0,1]$. All these three maps display chaotic behavior and are topologically conjugated to each other, thus sharing the same entropic spectrum. For a homogeneous partition of the interval with $p=2$ symbols $[0,1]=[0,1/2)\cup [1/2,1]$, the symbolic sequences of these maps don't have forbidden blocks, meaning that $h(0)=\log 2$. It is also well known that this partition is generating (not only that, it is also a Markov partition) and thus the entropies reach their supremum over all possible partitions. Furthermore, for these maps all symbol blocks are equiprobable, hence $h(\beta)=\log 2, \ \forall \beta$. In some sense, these maps are fully chaotic as their $p=2$ symbolic representation has strong randomness properties: one cannot distinguish these from a purely random sequence of two symbols. As a matter of fact for the binary shift there are actually no differences at all: the symbolic sequence in these cases not only induces a topological Markov chain but also a standard Markov chain, meaning that the sequence has no memory on previous states: the probability of a new symbol is $1/2$, irrespective of the previous history.\\

\noindent Now, the transition sequence has $p=4$ symbols, hence for the sake of comparison we firstly consider an homogeneous partition of $p=4$ symbols for the fully chaotic logistic map. Interestingly, in that case we find that this symbolic sequence has an entropic spectrum $h(\beta)\to \log 2 < \log 4, \ \forall \beta \geq 0$ as well. That means that we find forbidden patterns in this configuration.
Note that this is not unexpected and can be easily proved as follows: since the partition with two symbols is generating, the block entropy reaches its supremum over all partitions for that case, which is indeed $\log 2$, thus any other partition should yield entropies smaller or equal than $\log 2$. On the other hand, an homogeneous partition with $p=4$ symbols is a {\it refinement} of the partition with $p=2$ symbols, therefore invoking the monotonicity property of entropies over different partitions \cite{jost} we have that any refinement of the $p=2$ partition should yield an entropy larger or equal to the one for $p=2$, namely $\log 2$. Both conditions are only mutually satisfied when the inequality saturates, what concludes the proof.\\
The dynamical origin of the forbidden patterns in the fully chaotic logistic map for $p=4$ symbols can be explained in terms of the shape of the unimodal map $F(x)=4x(1-x)$. We will focus on a particular example to elucidate this relation. Define the symbols as the subintervals $A = [0,1/4)$, $B=[1/4,1/2)$, $C=[1/2,3/4]$ and $D=(3/4,1]$. Then for $m=2$, the pattern $(A,D)$ is forbidden (there are a total of 8 forbidden blocks of size 2). This means that in the symbolic sequence of the logistic map with $p=4$ symbols, the symbol $A$ will never precede the symbol $D$. The reason is simple: as $F$ is continuous it maps intervals into intervals, and together with the fact that $F(0)=0, \ F(1/4)=3/4$, this means that that $F:[0,1/4]\to [0,3/4]$, so starting from $A$ we will never reach $D$.
Interestingly enough, one can enumerate in this case the forbidden patterns and, surprisingly, we find that for every $m$ the amount of admissible and forbidden patterns exactly matches the formulas found for the type II null model. This suggests that one can make an enumeration of the forbidden patterns of the logistic map via a branching process similar to the one used to enumerate these patterns in the type II random process.\\

\noindent Once these synthetic processes have been detailed, let us consider the results obtained for the numerical experiments.

\begin{figure}[ht!]
\centering
\includegraphics[scale=0.37]{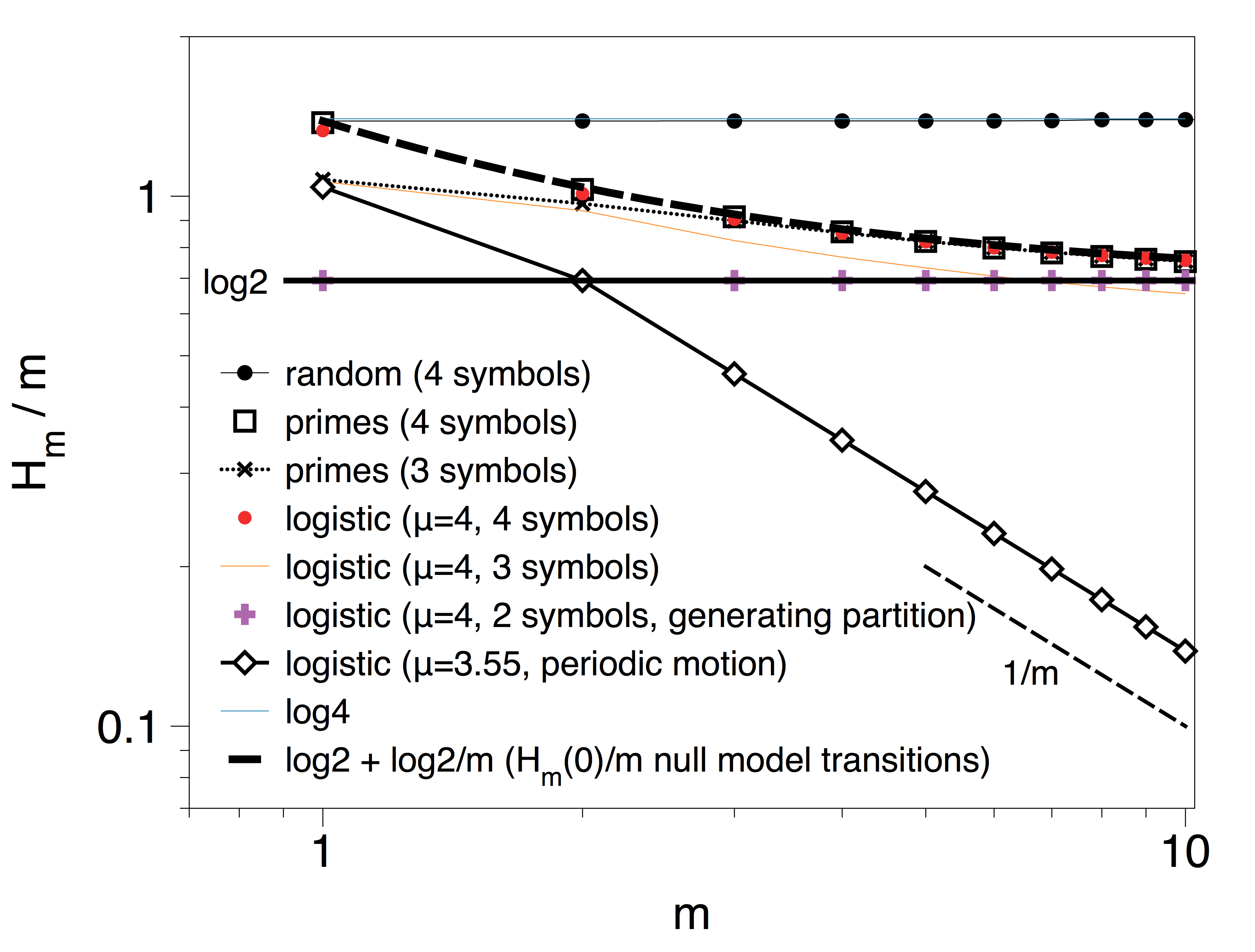}
\includegraphics[scale=0.37]{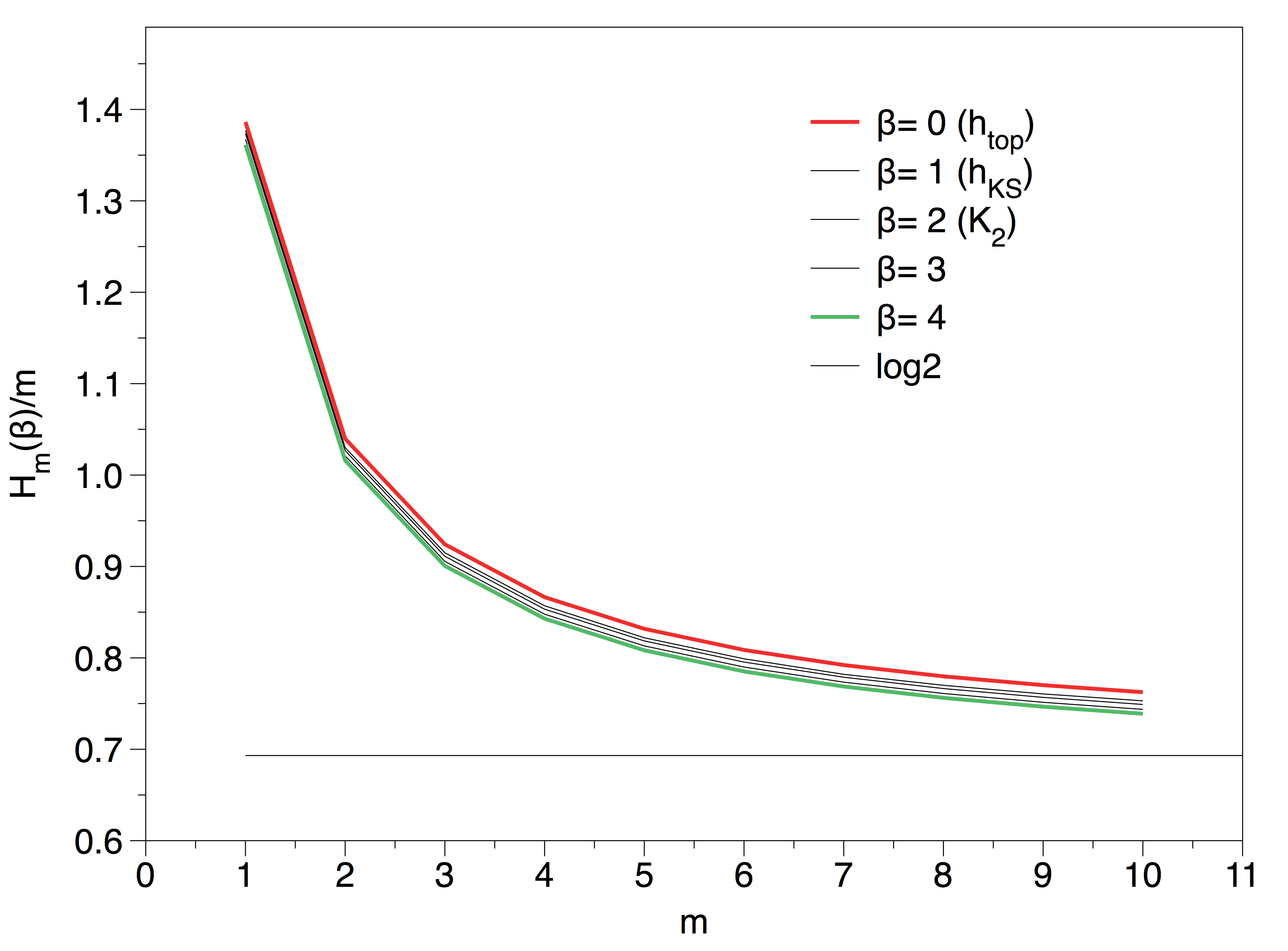}
\caption{(Left) Block entropies $H_m(1)/m$ for different symbolic sequences. (Right) Renyi block entropies $H_m(\beta)/m$ for different values of $m$ and $\beta$. All curves seem to converge to $\log 2$.}
\label{fig:h1}
\end{figure}

\subsection{Results and interpretation for the transition sequence}
Let us start by analysing the transition sequence. In the left panel of figure \ref{fig:h1} we show that, for the transition sequence with $p=4$ symbols, $H_m(1)/m \to \log 2$. This tendency is quantitatively similar for other values of $\beta$, see the right panel of the same figure, suggesting that $h(\beta)=\log 2$. This means that in the transition sequence there are forbidden patterns but those admissible occur equally often. This was indeed expected from the type II null model analysis displayed above. In the same figure we have plotted the analytical prediction for $H_m(0)/m=\log 2 + (\log 2)/m$ from the type II null model, which perfectly matches the numerical experiment. We conclude that the transition sequence shows strong randomness properties and cannot be distinguished from a purely random process.\\

\begin{figure}[ht!]
\centering
\includegraphics[scale=0.5]{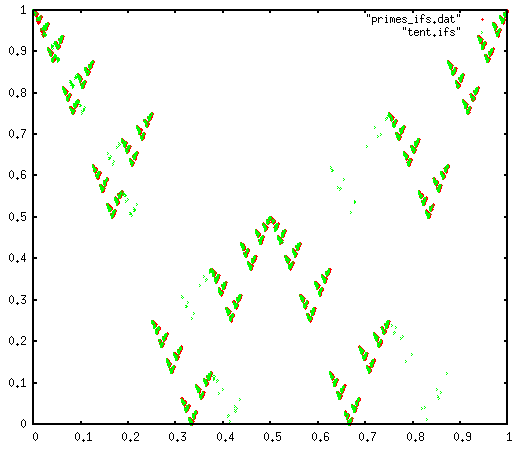}
\caption{IFS chaos game-like attractor for the prime transition sequence. The attractor has a fractal shape instead of being space-filling,  however a type II null model can account for such shape, as shown in figure \ref{fig:nullmodeltransition}.}
\label{fig:transitionIFS}
\end{figure}

Applying the chaos game (contraction factor $1/2$) to the transition sequence (figure \ref{fig:transitionIFS}) shows a clearly fractal shape. Interestingly, this is the same pattern that one finds for the symbolic sequence (4 symbols) generated by the tent map. As expected, this is different from a null model of type I but identical to a null model of type II (see right panel of figure \ref{fig:nullmodeltransition}), on agreement with the conclusion extracted from our entropic analysis.

\subsection{Results and interpretation for the Primes mod $k$}
\begin{figure}[ht!]
\centering
\includegraphics[scale=0.25]{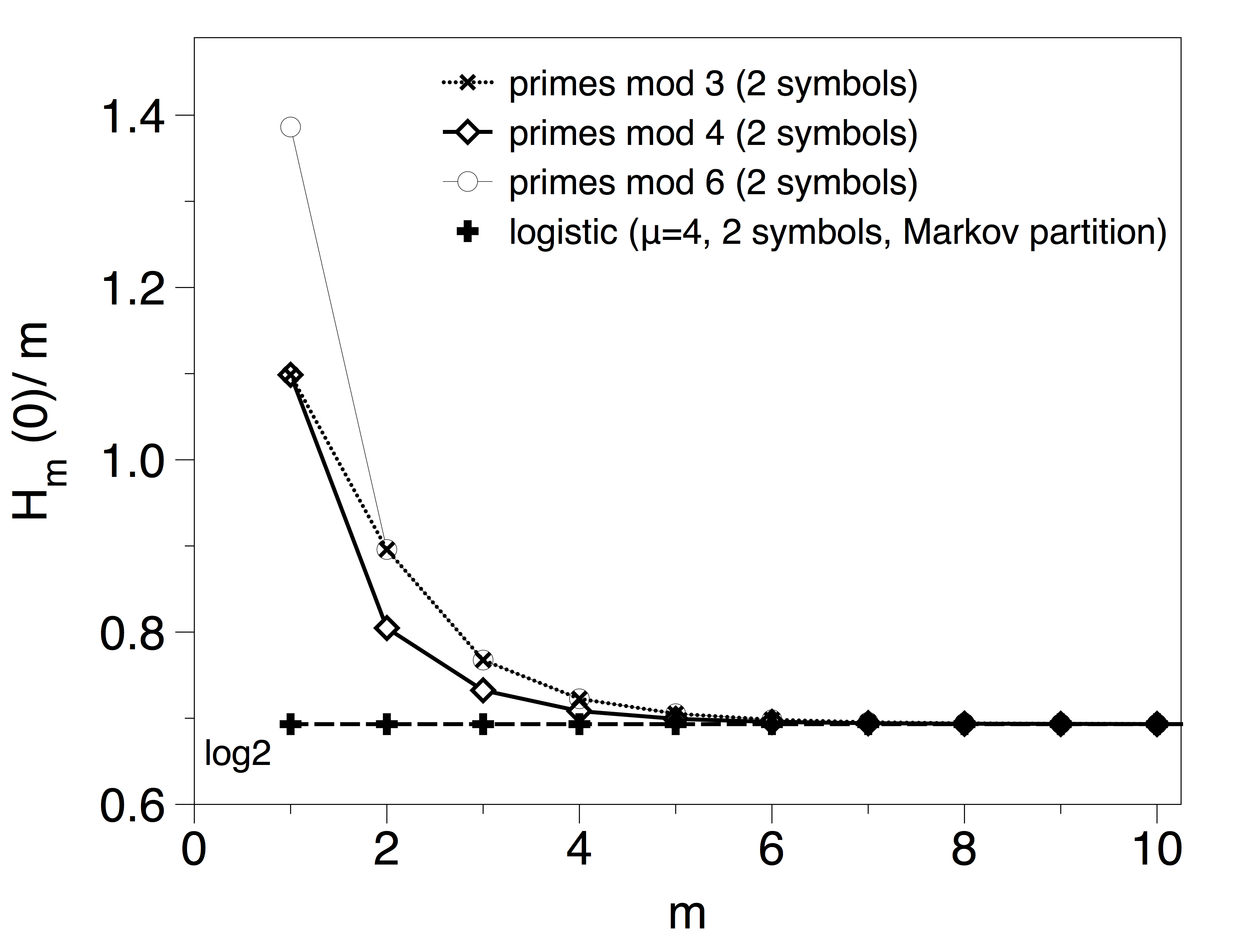}
\includegraphics[scale=0.25]{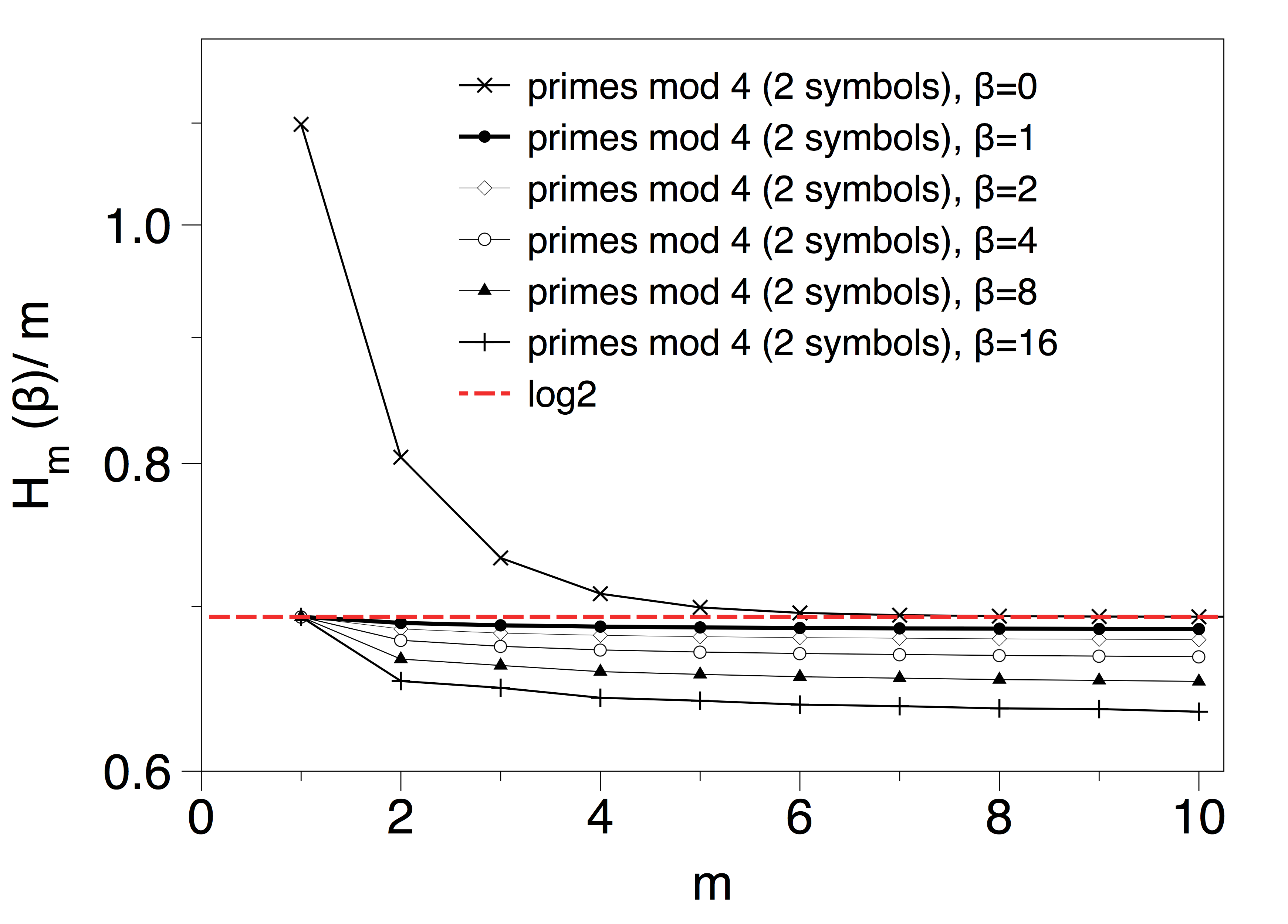}
\includegraphics[scale=0.25]{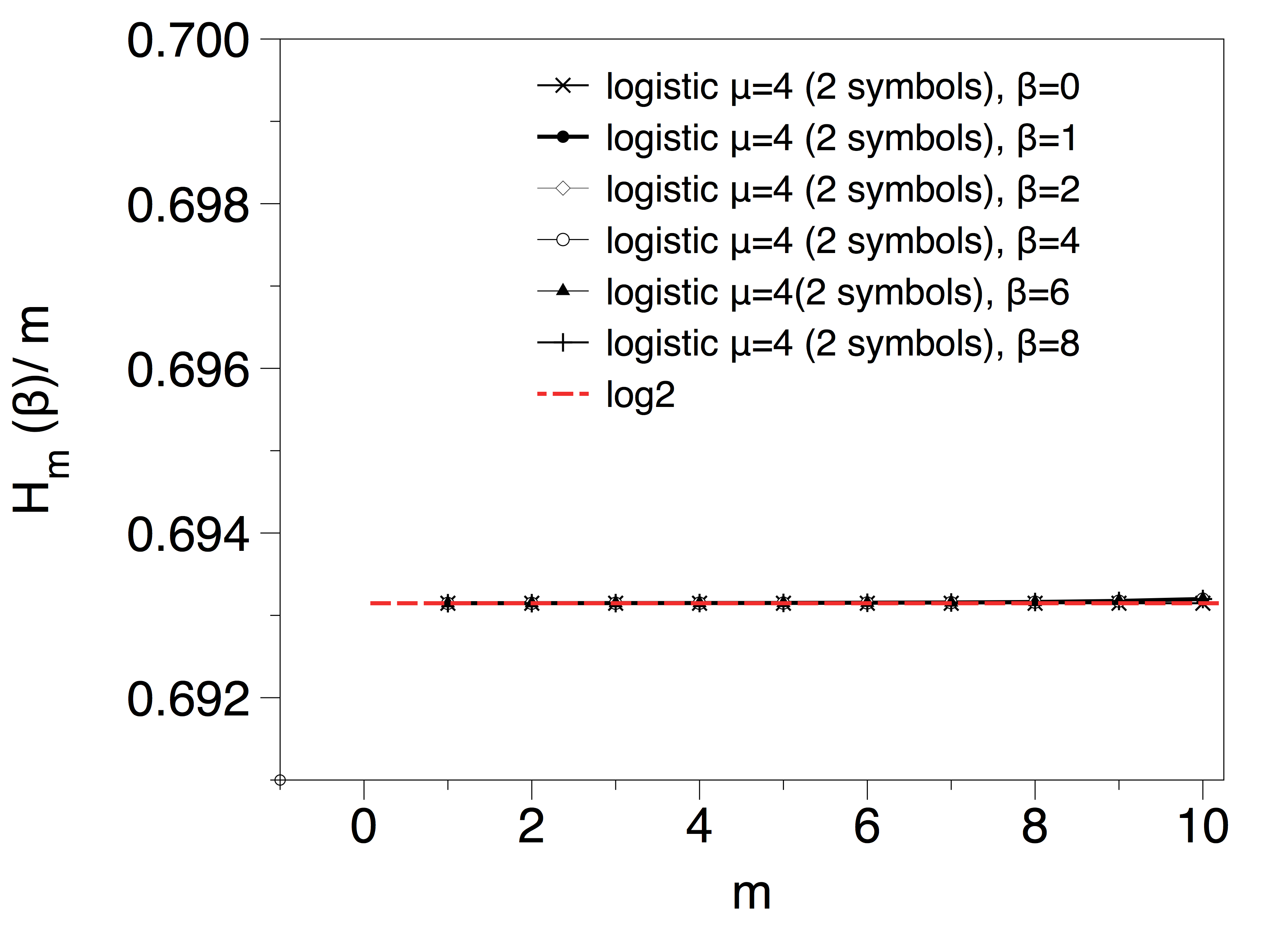}
\caption{(Left) Block entropies $H_m(0)/m$ for symbolic sequences primes mod k, with different values of $k$. We find that the topological entropy is always $\log 2$, suggesting that there are no forbidden patterns in primes modulo $k$. (Middle) Empirical Renyi block entropies $H_m(\beta)/m$ for different values of $m$ and $\beta$, extracted from the primes modulo 4. The spectrum does not collapse into a single value.(Right) Empirical Renyi block entropies $H_m(\beta)/m$ for different values of $m$ and $\beta$, extracted from the fully chaotic logistic map $x_{t+1}=4x_t(1-x_t)$ after symbolization with $p=2$ symbols. Results have been computed over a symbolic sequence of $N=10^6$ points, the same size as the one for the primes in the middle and left panels. The spectrum in this case collapses $h(\beta)=\log 2, \ \forall \beta$ as expected according to the topological conjugacy of the fully chaotic logistic map with the binary shift.}
\label{fig:h1_different_mods}
\end{figure}

The previous results suggest that the apparent lack of total randomness found in the transitions between Pythagorean and Gaussian primes was only due to the way the sequence was constructed, as a null model of type II found the same quantitative results. In order to gain a better understanding, we now consider the symbolic sequences generated via congruences of primes modulo $k$, where $k$ is such that $\phi(k)=2$, and we start by computing the spectrum of Renyi entropies $h(\beta)$ for different $k$ and different values of $\beta$, results are summarized in figure \ref{fig:h1_different_mods}.
In the left panel we consider the topological entropy for $k=3,4$ and $6$, for which there are two residue classes with infinite primes (two symbols). In every case we find a smooth convergence with the block size $m$ to an asymptotic value $h(0)\approx \log 2$, suggesting that there are no forbidden patterns in these sequences as expected. Incidentally, the reason why $H_1(0)=\log 3$ for $k=3,4$ and $\log 4$ for $k=6$, and that $\log 2$ is only reached asymptotically in every case is because depending on $k$, there are a few initial primes that do not fall into the residue classes with positive density: for $k=4$, $2\equiv 2 (\text{mod}\ 4)\neq 1,3$; analogously for $k=3$, $3 \equiv 0 (\text{mod}\ 3)\neq 1,2$; and finally for $k=6$ as this is a primorial with two factors we have two primes off the main residue classes $2\equiv 2 (\text{mod}\ 6),3\equiv 3 (\text{mod}\ 6) \neq 1,5$).\\

\noindent On the other hand, in the middle panel of the same figure we explore the spectrum of Renyi entropies for a particular modulo $k=4$ (results are analogous for $k= 3, 6$). As the set of transient symbols is finite, in the case of the KS and the rest of Renyi entropies that take into account the frequencies of each block this effect is smeared out.
The first important observation is that for $\beta=1$, the KS entropy converges to a value $0.685$ slightly below $\log 2\approx 0.693$. This small shift could be due to the fact that the density of $A$'s and $B$'s, while being asymptotically identical $p(A)=p(B)=1/2$ (Dirichlet), is different almost surely for any finite length sequence: again the famous Tchebytchev bias. Accordingly, a null model for this case does not yield an entropy equal to $\log 2$ but $-p(A)\log p(A)-p(B)\log p(B)$. To check whether this effect is simply a result of the Tchebychev bias, we have also run a type I null model where each symbol is drawn non-uniformly according to the empirical estimations of $p(A)$ and $p(B)$ for our experimental series. We find $p(A)\approx 0.4998$, $p(B)\approx 0.5002$ and $h(1)\approx 0.693$, which is too close to $\log 2$ to give account for the deviation observed for the KS entropy. Using a Cramer null model (type III) we also find the same result $\forall \beta, h(\beta)\approx 0.69$. The only possible solution is that, while there are no forbidden patterns and all blocks of size $m$ are admissible, they don't appear equally often and as $m$ increases this heterogeneity builds up stronger.\\
If that was actually happening, then we should expect that this non-uniformity gets more pronounced for higher values of $\beta$. As a matter of fact, and at odds with what we find for the binary shift map and the null models, we indeed find that $h(\beta)$ for the prime residues seems to be a monotonically decreasing function on $\beta$: for higher order values of $m$ even if all blocks are admissible, they are not equiprobable. To check that this is not a finite size effect, we have reproduced the same computation on a symbolic sequence of the same size extracted from the fully chaotic logistic map with $p=2$ symbols -which is known to yield $h(\beta)=\log 2$ when finite size effects are absent-. The spectrum in this case collapses to $h(\beta)=\log 2, \ \forall \beta$ as expected according to the topological conjugacy of the fully chaotic logistic map with the binary shift, and no relevant finite size effects are appreciated. Note that the Tchebytchev bias alone is not able as well to explain this systematic decrease (a non-uniform null model of type I gives $h(\beta)\approx 0.693 \ \forall \beta$). This confirms that in the case of the prime residues modulo $k$, we indeed have a nontrivial spectrum of Renyi entropies due to the different rates of appearance of every block of size $m$. For $m=1$ this systematic non-uniformity is precluded by virtue of the Dirichlet theorem, that assigns the same density for all admissible residue classes, so this is a higher order pattern.\\
We use a proxy to $h(\beta)\approx H_{10}(\beta)/10$ and we plot this value as a function of $\beta$ in figure. We deduce that metric does indeed play a substantial role, at odds of what was perceived in the $\{AA,AB,BA,BB\}$ symbolization. The reason for the uneven appearance of each block will be elucidated in the next section, where we consider the gaps residues.

\begin{figure}[ht!]
\centering
\includegraphics[scale=0.45]{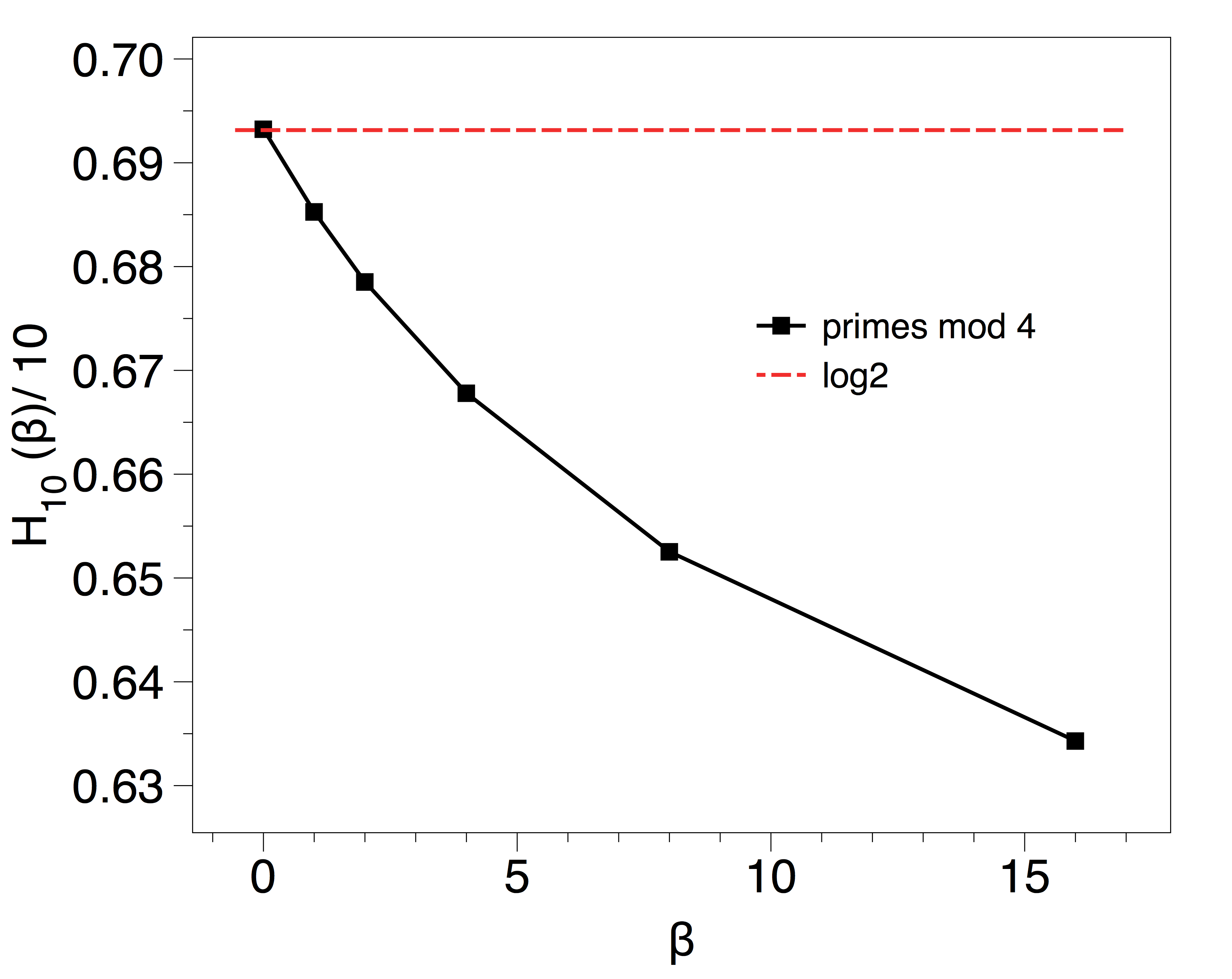}
\caption{Spectrum of Renyi entropies associated to the primes modulo 4.}
\label{fig:spectrum}
\end{figure}

\subsection{Results and interpretation for the gaps residue sequence}
Here we focus on the sequence of prime gaps, with $\text{gap}(n)=p(n+1)-p(n)$. In some sense gaps are the derivative of the primes. As it happens in complex chaotic systems such as turbulent fluids, it is not position but velocity (the derivative field) the observable which shows a richer structure.\\
Gaps are always even numbers and as there exist arbitrarily large gaps, in order to make the gap sequence stationary we consider this sequence modulo $6$. Remember that the choice of the modulus is not arbitrary here, as the residue classes classify gaps into three big families: 2 mod 6 (twin-like), 4 mod 6 (cousin-like) and 0 mod 6 (sexy-like).
\subsubsection{Entropic analysis}
\begin{figure}[h]
\centering
\includegraphics[scale=0.37]{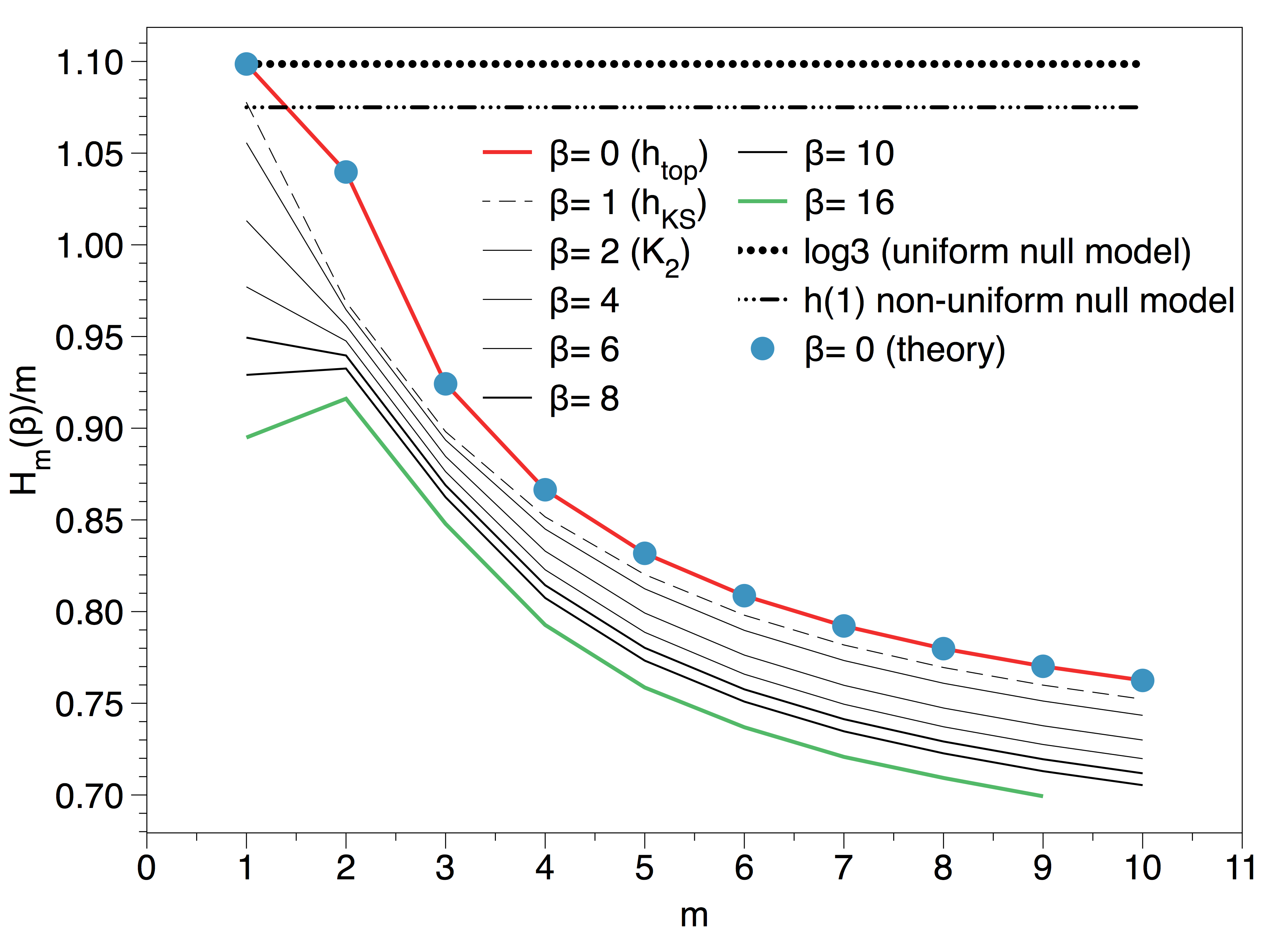}
\includegraphics[scale=0.37]{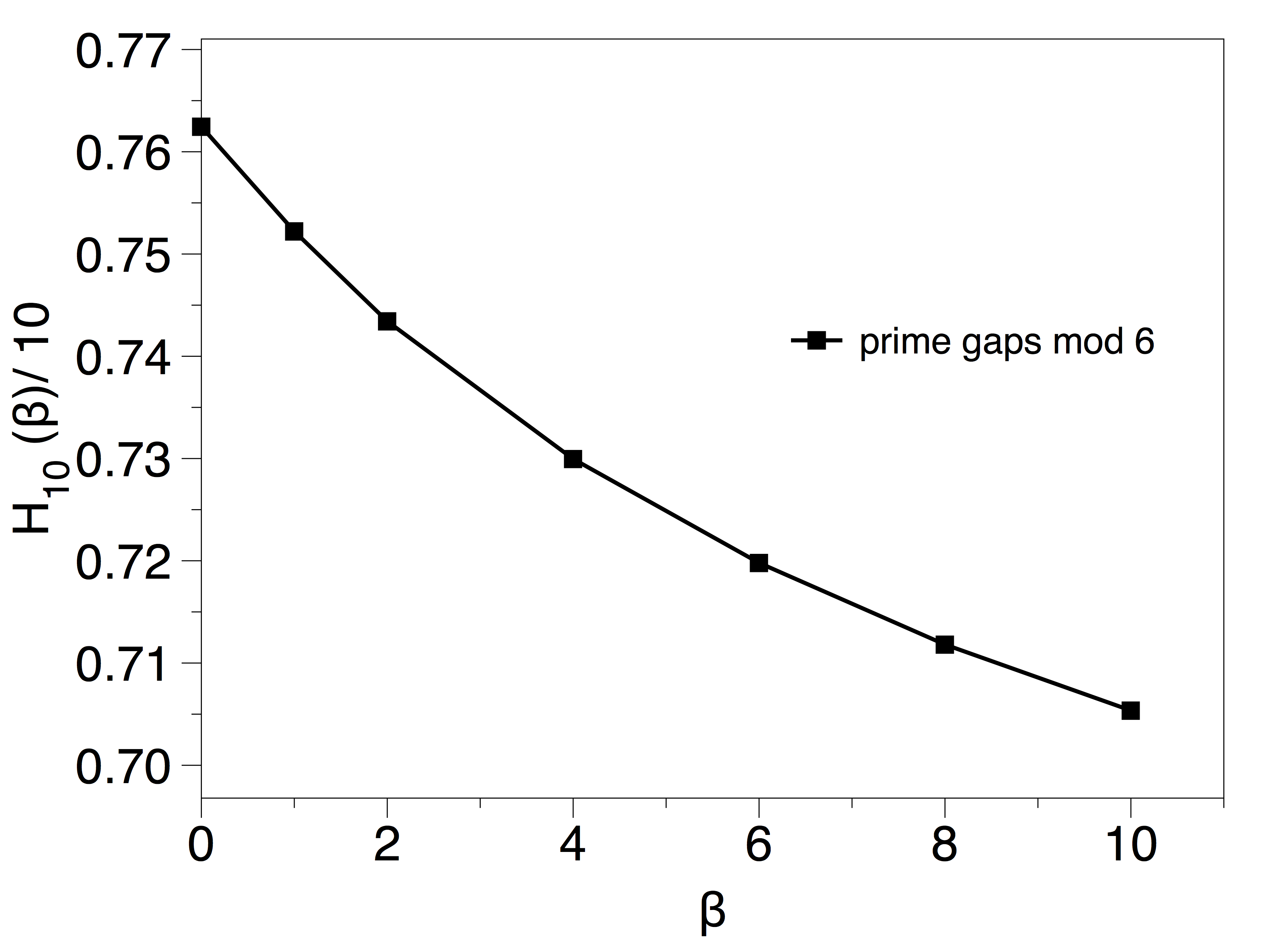}
\caption{(Left panel) $H_m(\beta)/m$ associated to the primes gaps modulo 6. (Right panel) Spectrum of Renyi entropies associated to the primes gaps modulo 6.}
\label{fig:spectrumgaps}
\end{figure}

In principle, a null model for this sequence should be of the type I, yielding for the uniform case $h(\beta)\approx \log 3$. For the more realistic case where the empirical frequencies of the residue classes are considered, we have observed that for the sequence of the first $10^6$ gaps,
\begin{equation}
p(0)\approx 0.43,\  p(2)\approx p(4)\approx 0.28.
\label{observed}
\end{equation}
For these frequencies, we find  $H_1(1)\approx 1.075$. The null model predicts $h(\beta)=H_1(\beta)$, so the first partial entropy should be enough to extrapolate the asymptotic value.
Using a Cramer null model we also find the same result. However, the actual spectrum is shown in figure \ref{fig:spectrumgaps}, pointing out two unexpected features: (i) the convergence of Renyi entropies to finite, positive values smaller than those predicted by the null models. This happens even for $\beta=0$, pointing out to the presence of forbidden patterns in this sequence: blocks of $m$ symbols that  appear in the null models but are forbidden in the prime gaps. And (ii) there seems to be a non-trivial, monotonic dependence on $\beta$ which is not found in the null models. In what follows we discuss both results.\\

\noindent {\bf Explaining forbidden patterns. }Some observed low order forbidden patterns are enumerated in table \ref{tab:forbiddengaps}. The first forbidden pattern is for a block of size $m=2$ and consists in $(4,4)$. A forbidden pattern in the gaps residue sequence such as $(4,4)$ actually relates to an infinite set of forbidden patterns in the prime sequence. First, all gap pairs $(f,f')$ with $f=6n+4, \ f'=6n'+4, \ \forall n,n'$ are congruent to $(4,4)$, that is, consecutive gaps such as $(4,4)$, $(4,10)$, $(4,16)$, $(10,4)$, etc are all forbidden. In the prime sequence, each of these forbidden gap pairs is associated with a forbidden prime triple of the form $(q,q+f,q+f')$, with $q$ prime.\\

\begin{table}[ht]
\centering
\begin{tabular}{|c|c|c|}
\hline
m & {\bf ${\cal F}(m)$ }  & $|{\cal F}(m)|$\\
\hline
\hline
1  &  $\emptyset$ & 0\\
\hline
2  & \{(4,4)\} & 1\\
\hline
3 & $\{(0,2,2),(0,4,4),(2,0,2),(2,2,0),(2,2,2),(2,4,4),(4,0,4),(4,2,2),(4,4,0),(4,4,2),(4,4,4)\}$& 11\\
\hline
4 & $\{(0,0,2,2),(0,0,4,4),(0,2,0,2),(0,2,2,0),(0,2,2,2),(0,2,2,4),(0,2,4,4),(0,4,0,4),(0,4,2,2),(0,4,4,0),$&49\\
&$(0,4,4,2),(0,4,4,4),(2,0,0,2),(2,0,2,0),(2,0,2,2),(2,0,2,4),(2,0,4,4),(2,2,0,0),(2,2,0,2),(2,2,0,4),$&\\
&$(2,2,2,0),(2,2,2,2),(2,2,2,4),(2,2,4,0),(2,2,4,4),(2,4,0,4),(2,4,2,2),(2,4,4,0),(2,4,4,2),(2,4,4,4),$&\\
&$(4,0,0,4),(4,0,2,2),(4,0,4,0),(4,0,4,2),(4,0,4,4),(4,2,0,2),(4,2,2,0),(4,2,2,2),(4,2,2,4),(4,2,4,4),$&\\
&$(4,4,\dots,\dots)\}$&\\
\hline
\end{tabular}
\caption{Set of forbidden blocks  ${\cal F}(m)=\{(s_1\dots s_m), s_i\in\{0,2,4\}\}$ of size $m=1, 2, 3, 4$ in the sequence of gap residues modulo $6$.}
\label{tab:forbiddengaps}
\end{table}
\noindent From a dynamical systems viewpoint, these results suggest that the gap residue sequence manifests chaotic behavior (as the amount of admissible blocks grows exponentially with $m$ and thus the KS entropy is positive), but of weaker intensity as $h(0)<\log p$. According to Pesin identity \cite{beck}, for symbolic dynamics with positive KS entropy, this quantity is equal to the Lyapunov exponent of the underlying dynamical system, this latter quantity being responsible for the exponentially fast separation of initially close trajectories. In a system with $p$ symbols, information could be erased at least as fast as $\log p$ per time step (in Information Theory both the topological and KS entropies are usually defined in base 2 rather than in base $e$, and in this setting information can be erased as fast as $\log_2 2= 1$ bit per iteration). In any case, for the gaps residue sequence modulo 6 we have $p=3$, so
the dynamics is chaotic albeit information is lost at a slower pace as $h(0)\to \log 2< \log 3$.\\

\noindent Why is this the case? Why is the dynamics underlying the prime gap residue sequence 'less chaotic' than a shift of finite type or a purely random process?
In other words: what is the nature of these forbidden blocks? Consider for instance the block $(2,2)$, which is not forbidden as it appears at least once in the prime progression $3,5,7$. It is easy to show that such progression is actually the only possibility. The proof consists in studying the divisibility properties of $q,q+2$ and $q+4$. It turns out that in this progression there is always an element which is not prime because it is a multiple of three. Assume for a contradiction that this is not the case, and consider the remainder of the integer division $q \mod3$. If $q$ is not a multiple of three, the remainder should either be 1 or 2. If the remainder is 1, then $q+2$ is a multiple of three. If the remainder is 2, then $q+4$ is a multiple of three. The only case where there is no contradiction is when $q=3$ (a prime which is multiple of three), what finishes the proof.\\
Note that this proof also certifies that no progression $q\ \mod 6,\ (q+2) \ \text{mod}\ 6,\ (q+4) \ \text{mod}\ 6$ other than 3, 5, 7 is possible, i.e. this holds not just for a twin triplet but for any twin-like triplet of the type $q,q+6n +2, q+6n'+4$ where $n'\geq n$.\\
The origin of forbidden blocks can be thus directly linked to the divisibility properties of the integers. Actually, it is well-known that a block of $m$ consecutive gaps $(2g_1,\dots,2g_m)$, which gives rise to a prime block of the type $q,q+2g_1,q+2g_1+2g_2,\dots,q+\sum_{i=1}^m 2g_i$ is forbidden if and only if one can find a prime $r$ for which each all and every partial series $\sum_{i=1}^{t\leq m}  g_i$ is congruent to a different residue from $1,2,\dots, r-1$. For instance, in the case above of a gap duple $(2,2)$ that generates a prime sequence of the form $q,q+2,q+4$, $2\equiv 2 \mod 3$ and $4\equiv 1 \mod 3$.
For $(4,4)$ again all residues for $r=3$ are ticked, thus $(4,4)$ is forbidden. Using this principle one can therefore systematically enumerate these forbidden patterns. For a generic $m$ one finds
$$|{\cal A}(m)|=2^{m+1};\ |{\cal F}(m)|=3^m-2^{m+1}$$
from which we derive the partial entropies $H_m(0)$ for the topological entropy
\begin{equation}
\frac{H_m(0)}{m}=\frac{1}{m}\log |{\cal A}(m)|=\Big (1+\frac{1}{m}\Big )\log 2
\end{equation}
for $m\geq 2$ and $H_1(0)=\log 3$. This monotonic decay is in good agreement with the experiment (see the left panel of figure \ref{fig:spectrumgaps}). Interestingly, the number of admissible blocks of size $m$ -which is here related to the divisibility properties of the integers- precisely coincides with the number of admissible blocks in the type II null model, see eq. \ref{nullII}.
Altogether,
$$h(0)=\lim_{m\to \infty}\frac{H_m(0)}{m}=\log 2$$

\noindent {\bf Monotonic dependence on $\beta$ and the distribution of blocks.} In order to find an analytical expression for $h(\beta>0)$ which would allow us to elucidate whether the Renyi spectrum is indeed non-trivial or, on the contrary, whether this is just a finite size effect which while not present in the null models might be present for finite size statistics of gaps residue sequences but vanish asymptotically. We would need to be able to find an analytical expression for the frequencies of each admissible block.\\
Let us start by considering a Cramer null model. According to this model any large integer $q$ has roughly a probability $1/\log q$ of being a prime, then assuming probability independence the probability that an $m$-tuple of integers is prime is simply $1/(\log q)^{m+1}$. However, obviously the Cramer model predicts the same probability for every m-tuple, that is, in a Cramer random model every gap block of size $m$ would be equiprobable and therefore the frequency of each admissible block of size $m$ is simply $1/|{\cal A}(m)|$. Under this situation it is easy to show that $H_m(\beta)=H_m(0)=\log |{\cal A}(m)|$: a Cramer null model does not explain the dependence on $\beta$ found in the gap residue sequence.\\
Fortunately however, a well-known conjecture in number theory comes to save the day. First, let us define a prime $m$-tuple as a sequence of consecutive primes of the form $p,p+2g_1,p+2g_2,\dots,p+2g_m$. Such a prime $m$-tuple is trivially associated to a gap block $(2g_1,\dots,2g_m)$ (of course, in our case many different gap blocks have the same associated residue block). The diameter of a prime m-tuple is the difference of its largest and smallest element, i.e. $2g_m$. For a fixed $m$, the $m$-tuple with smallest diameter is called a prime constellation. For instance, for a prime duple $(p,p+2g_1)$ one can find twin primes ($g_1=1$), cousin primes ($g_1=2$), sexy primes ($g_1=3$) and so on, all associated to a gap block $(2g_1)$. The smallest diameter is for $g_1=1$ and therefore for $m=1$ the only prime constellation consist of twin primes.
For gap blocks of size $m=2$, the smaller diameter is $2g_2=6$, and there are two possible constellations for that, associated with the gap blocks $(2,4)\ (\equiv (2,4) \mod 6)$ and $(4,2)\ (\equiv (4,2) \mod 6)$. An example of the former is the prime triple $(5, 7, 11)$ whereas for the latter the smallest constellation is $(7,11, 13)$.\\
It is clear that prime constellations generate only a subset of our gap residue blocks, but still an infinite subset, and more importantly, for a fixed $m$ there is in general more than one constellation. Now, the celebrated m-tuple conjecture by Hardy and Littlewood \cite{hardy} states that the frequencies of these prime constellations can be computed explicitly:\\

\noindent {\bf Prime $m$-tuple conjecture}.\\
The amount of prime constellations $[p, p+2g_1, \dots, p+2g_m]$ found for $p\leq x$ is given asymptotically by
$$\pi_{g_1,g_2,\dots,g_m}(x)\sim C(g_1,\dots,g_m)\int_2^x\frac{dt}{(\log t)^{m+1}},$$
where
$$C(g_1,\dots,g_m)=2^m\prod_{q>2} \frac{1-\omega(q;g_1,\dots,g_m)/q}{(1-1/q)^{m+1}}$$ are the so-called Hardy-Littlewood constants, the product runs over odd primes $q>2$ and the function $\omega(q;g_1,\dots,g_m)$ denotes the number of distinct residues in $0,g_1,\dots,g_m$ mod $q$.

\noindent Incidentally, note that for readability we have used the symbol $m$ to be consistent with the previous exposition, however this conjecture is usually stated as the $k$-tuple conjecture. Remarkably, these probabilities are proportional to the Cramer model but, contrary to the Cramer model, do indeed depend on the particular block type via the Hardy-Littlewood constants.
This enables us to explore the uneven distribution of blocks via this conjecture:\\

\noindent For $m=1$, we are considering pairs of primes. The prime constellation in this case, as we know, is for $g_1=1$, therefore the conjecture predicts that the amount of twin primes below $x$ is
$$\pi_2(x)\sim C(1)\int_2^x\frac{dt}{(\log t)^{2}},$$
where $$C(1)=2\prod_{q\ \text{prime}\ > 2} \frac{1-2/q}{(1-1/q)^2}\approx 1.320324$$
The formula given by Hardy and Littlewood can actually be used to estimate the amount of any prime $m$-tuples, not necessarily only those with smaller diameter. For instance, still for $m=1$ we can work out the case for which $g_1=2$ (cousin primes) and $g_1=3$ (sexy primes), finding
$$C(2)=C(1), \ C(3)= 2C(2), \ C(g_1)=C(1)\cdot\prod_{q|g_1}\frac{q-1}{q-2}.$$
From this first analysis it is straightforward that sexy pairs are, asymptotically, twice more frequent than twins and cousins, what already suggests for $m=1$ that we have a non-uniform distribution of blocks.\\

\noindent The uneven distribution of blocks can be further assessed by enumerating for a given size $m$ all the $|{\cal A}(m)|$ Hardy-Littlewood constants $\{C(g_1,\dots,g_m)\}$. This is in principle possible but is a formidable exercise in practice. An additional technical issue to have in mind is that here we deal with blocks of prime gap residues, so for a given $m$ we need to sum over all indices congruent with a given residue block. Let us consider again the case $m=1$. The amount of gaps congruent to 0 mod 6 (sexy primes) is given by taking into account gaps 6, 12, 18, etc so $g_1=3,6,9,$ etc. Accordingly, this amount is associate to $\sum_{i=1}^{\infty} C(3i)$. Formally, if we define the normalization factor ${\cal Z}$ as: 

$${\cal Z}=\sum_{i=1}^{\infty}C(3i) + \sum_{i=0}^{\infty} C(3i+1) + \sum_{i=0}^{\infty} C(3i+2)$$

the densities of each gap residues is therefore:
$$p(0)=\frac{1}{{\cal Z}}\sum_{i=1}^{\infty} C(3i), \ p(2)=\frac{1}{{\cal Z}}\sum_{i=0}^{\infty} C(3i+1), \ p(4)=\frac{1}{{\cal Z}} \sum_{i=0}^{\infty} C(3i+2)$$
Truncating these series at order $i$ accounts for the gaps of size at most $2i$. For instance, for third order (accounting for gaps up to size 18) we have an approximation:
$$p(0)\approx \frac{C(3)+C(6)+C(9)}{\sum_{i=1}^9 C(i)}\approx 0.479, \
p(2)\approx \frac{C(1)+C(4)+C(7)}{\sum_{i=1}^9 C(i)}\approx 0.255, \
p(4)\approx \frac{C(2)+C(5)+C(8)}{\sum_{i=1}^9 C(i)}\approx 0.266,$$
which is still reasonable far away from the empirical frequencies (eq. \ref{observed}) (at order four $p(0)\approx 0.471$ so convergence is slow). Still, this analytical approximation already certifies that the three possible gap residue blocks of size $m=1$ are not uniformly represented. Formally, one can express, assuming Hardy-Littlewood conjecture, the asymptotic probability of an admissible gap residue block of size $m$ $(2g_1,2g_2,\dots,2g_m)\in {\cal A}(m)$, where $g_i\in \{0,1,2\}$ as
\begin{equation}
p(g_1,g_2,\dots,g_m)=\frac{1}{{\cal Z}}\sum_{i_1}\sum_{i_2}\cdots \sum_{i_m} C(3i_1 + g_1, 3i_2 + g_2, \dots, 3i_m + g_m)
\label{general}
\end{equation}
where the summation in $i_j$ runs between $i_{j-1}$ and $\infty$ when $g_j=1,2$ and $g_{j-1}=0$ and  between $1+i_{j-1}$ and $\infty$ otherwise, and ${\cal Z}$ is the normalization factor. Eq. \ref{general}, together with eqs. \ref{ks} and \ref{renyi} constitute the formal solution to the problem.\\

\noindent Let us consider the case $m=2$, for which there are 8 admissible blocks. To prove non-uniformity for $m=2$ we only need to find that two different blocks have different frequencies: for convenience we will concentrate on the first two ones, namely $(2g_1,2g_2)=(0,0)$ in $\mod 6$ (which gathers all prime triples whose gaps are a multiple of 6) and $(2g_1,2g_2)=(0,2)$. Eq. \ref{general} reduces in these cases to:\\
$$p(0,0)=\frac{1}{{\cal Z}}\sum_{i_1=1}^{\infty}\sum_{i_2=i_1+1}^{\infty} C(3i_1, 3i_2)\ \textrm{and} \ p(0,2)=\frac{1}{{\cal Z}}\sum_{i_1=1}^{\infty}\sum_{i_2=i_1}^{\infty} C(3i_1, 3i_2+1)$$
where
\begin{eqnarray}
&&{\cal Z}=\sum_{i_1=1}^{\infty}\sum_{i_2=i_1+1}^{\infty} C(3i_1, 3i_2) + \sum_{i_1=1}^{\infty}\sum_{i_2=i_1}^{\infty} C(3i_1, 3i_2+1)
+\sum_{i_1=1}^{\infty}\sum_{i_2=i_1}^{\infty} C(3i_1, 3i_2+2)+\sum_{i_1=0}^{\infty}\sum_{i_2=i_1+1}^{\infty} C(3i_1+1, 3i_2)\nonumber \\
&&+\sum_{i_1=0}^{\infty}\sum_{i_2=i_1+1}^{\infty} C(3i_1+1, 3i_2+1)
+\sum_{i_1=0}^{\infty}\sum_{i_2=i_1+1}^{\infty} C(3i_1+1, 3i_2+2)
+\sum_{i_1=0}^{\infty}\sum_{i_2=i_1+1}^{\infty} C(3i_1+2, 3i_2)\nonumber \\
&&+\sum_{i_1=0}^{\infty}\sum_{i_2=i_1+1}^{\infty} C(3i_1+2, 3i_2+1) \nonumber
\end{eqnarray}
 is the normalization sum.
At second order in $i_j$,
$$p(0,0)\approx\frac{C(3,6)+C(3,9) +C(6,9)+C(6,12)}{{\cal Z}|_2},$$
and
$$p(0,2)\approx\frac{C(3,4)+C(3,7) +C(6,7)+C(6,10)}{{\cal Z}|_2},$$
where ${\cal Z}|_2$ is the normalization sum truncated at order two. At this point we can make use of Hardy-Littlewood's conjecture one more time. If we label
$$a:= \prod_{q\ \text{prime} \ \geq 5} \frac{1-2/q}{(1-1/q)^3}$$
and
$$b:= \prod_{q \ \text{prime} \ \geq 5} \frac{1-3/q}{(1-1/q)^3},$$
then after a lengthy but trivial computation we get $C(3,6)=9b$, $C(3,9)=9b$, $C(6,9)=9b$, $C(6,12)=9b$ hence
$$p(0,0) \approx \frac{1}{{\cal Z}|_2}36b$$
whereas $C(3,4)=9b/2$, $C(3,7)=45b/8$, $C(6,7)=45b/8$, $C(6,10)=27b/4$ thus
$$p(0,2)\approx \frac{1}{{\cal Z}|_2}\frac{1}{2}45b \neq p(0,0),$$
which is enough to conclude non-uniformity in the approximation of the frequencies of blocks of size $m=2$. This is a clear support to the apparence monotonic dependence on $\beta$ of the Renyi entropies. A full proof would require to perform the infinite summations in each case, what appears to be a quite challenging endeavor and is left as an open problem.
\subsubsection{IFS}
\begin{figure}[ht!]
\centering
\includegraphics[scale=0.4]{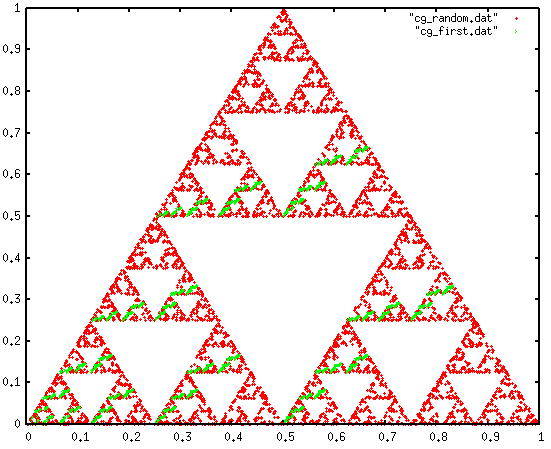}
\includegraphics[scale=0.25]{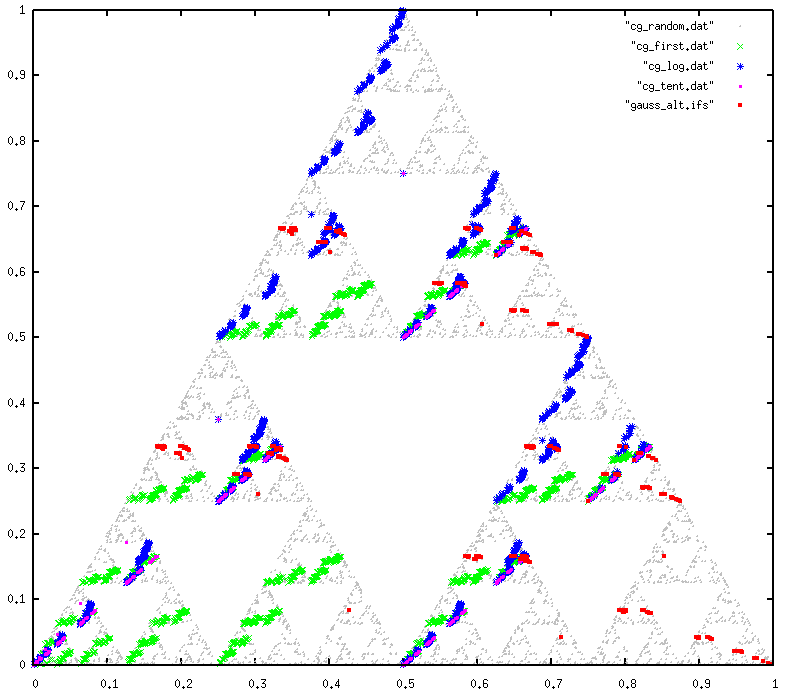}
\caption{(Left panel) IFS chaos game-like attractor for the type I null model with $p=3$ symbols. In this case the attractor is the Sierpinski triangle. (Right panel) IFS chaos game-like attractor for the prime gap sequence and a symbolic sequence extracted from other chaotic maps. The IFS associated to the prime gap sequence is a subset of the attractor. The similarity with the tent map is again notable, where in this case the attractor of the  tent map is a subset of the attractor of the prime gap sequence. As a comparison, we have included a different chaotic map (the Gauss map) which shows no similarities with the gaps.}
\label{fig:IFSgaps}
\end{figure}

The analysis performed a la IFS is shown in figure \ref{fig:IFSgaps}, where we also plot the result of a type I null model. We find that the prime gaps selects a subset of the attractor generated by the null model. This particular subset could be put in correspondence with the subset of admissible sequences found above, and in some sense the IFS attractor plays the geometric role of the distribution of forbidden patterns.

\section{Discussion}
The cross-disciplinary transfer of concepts and tools has always been a fruitful strategy for finding unexpected and unorthodox ways of addressing problems, which has usually generated new knowledge. In this work we have tried to illustrate this idea by making using the concepts, focus and techniques used traditionally in nonlinear dynamics and complexity science to explore the structure of certain sequences appearing in number theory. We have considered three types of sequences that emanate from the prime number sequence: the transition between Pythagorean and Gaussian primes (transition sequence), the sequence of prime residues modulo $k$ (with $\phi(k)=2$) and the sequence of prime gap residues modulo 6 (inducing three symbols that can be associated with twin-like, cousin-like and sexy-like pairs). All these sequences are stationary and it is easy to prove using Dirichlet's theorem for arithmetic progressions that each symbol appears infinitely often, hence they are amenable to a symbolic dynamics analysis. We have thus considered these symbol sequences as if they were generated by a hidden dynamical system whose trajectories were symbolized via some unknown partition of the system's phase space and have consequently explored some dynamical properties of the system via a symbolic analysis of the sequences.\\

\noindent First, in every sequence we have found compelling evidence of chaotic behavior, as given by a positive KS entropy. In the case of the transition sequence, we have found that the sequence has a lower entropy than expected from a purely random process or a shift of finite type but we have shown that this decrease is due to the onset of spurious forbidden patterns, not associated to the underlying dynamics but to the way of defining transitions between primes. The IFS attractor generated following the transition sequence is a fractal that coincides with the equivalent attractor for the tent map and for a properly defined null model, concluding that this transition sequence has strong randomness properties, as expected.\\

\noindent As for the residues of the primes modulo $k$, we have found that this sequence is maximally chaotic in the sense that its topological entropy matches the analogous one found for the binary shift map. While lacking forbidden patterns and at odds with the behavior of the binary shift map, we have found that this sequence displays a non-trivial spectrum of Renyi entropies which unexpectedly suggest that every symbol block of size $m>1$, while admissible, occurs with different probability. This non-uniform distribution of blocks for $m>1$ contrasts Dirichlet's theorem that guarantees equiprobability for $m=1$.\\

\noindent Finally, we have  explored in a similar fashion the sequence of prime gap residues. This sequence is again chaotic (positivity of Kolmogorov-Sinai entropy), however chaos is weaker in this case as we have found non-spurious forbidden patterns for every block of size $m>1$ that yield entropies lower than those for a null model or a shift of finite type with the same number of symbols. We were able to relate the onset of these forbidden patterns with the divisibility properties of integers. Interestingly, the amount of admissible blocks is precisely the same, for a given block size, than both the type II null model and the fully chaotic logistic map with $p=4$ symbols. This suggests that a systematic enumeration of admissible sequences in the gap sequence could be figured out without resorting to any number-theoretic properties, something that is left as an open problem. We have also given a geometric visualization of this phenomenon by showing that the IFS attractor associated to this sequence is only a subset of the Sierpinski triangle, which is the attractor for a three-vertex Chaos Game.\\
Moreover, we have found again that this sequence displays a monotonically decreasing spectrum of Renyi entropies, which suggests that gap block residues of size $m>1$ are not equiprobable. A Cramer random model  cannot explain this apparent dependence on $\beta$ as for that model, given $m$ each block of size $m$ is by construction equiprobable. However,
the frequencies of these blocks can be computed explicitly assuming Hardy-Littlewood's $k$-tuple conjecture. We have built on this relation to give analytical arguments and approximations that support the monotonic dependence of the Renyi entropies on $\beta$, however an exact and closed form expression for these entropies has remained elusive and constitutes an interesting open problem.

%\bibliography{paperbib}{}
\bibliographystyle{apsamp}

\end{document}